\input ./style/arxiv-general.cfg
\documentclass[MSNbibl,number,citesort,dvips]{arxbj}
\makeatletter
   \@ifpackageloaded{graphicx}{}{\usepackage{graphicx}}
\makeatother
\usepackage{eufrak}
\usepackage{upgreek}
\usepackage{dcolumn}

\aid{0}
\volume{22}
\issue{1}
\pubyear{2016}
\firstpage{444}
\lastpage{493}
\doi{10.3150/14-BEJ665} 
\docsubty{FLA}

\makeatletter
\newcolumntype{d}[1]{D{.}{.}{#1}}
\newcommand{\rrvert}{\vert}
\newcommand{\llvert}{\vert}
\newcommand{\notag}{\nonumber}
\newcommand{\eqref}[1]{(\ref{#1})}
\newtheorem{lemma}{Lemma}
\newtheorem{proposition}{Proposition}
\newtheorem{theorem}{Theorem}
\newtheorem{corollary}{Corollary}

\newremark{remark}{Remark}
\newcommand{\cm}{{\EuFrak B}}

\newcommand{\ve}[1]{\mathbf{#1}}
\newcommand{\veee}[1]{\bolds{#1}}
\makeatother

\begin{document}
\begin{frontmatter}

\title{Tail behavior of sums and differences of log-normal random variables}
\runtitle{Tails of log-normal sums and differences}

\begin{aug}
\author[1]{\inits{A.}\fnms{Archil}~\snm{Gulisashvili}\thanksref{1}\ead[label=e1]{gulisash@ohio.edu}}%
\and
\author[2]{\inits{P.}\fnms{Peter}~\snm{Tankov}\corref{}\thanksref{2}\ead[label=e2]{tankov@math.univ-paris-diderot.fr}}
\dedicated{This article is dedicated to the memory of Peter
Laurence.}
\address[1]{Department of Mathematics, Ohio University, USA. \printead{e1}}
\address[2]{Laboratoire de Probabilit\'es et Mod\`eles
Al\'eatoires,
Universit\'e Paris Diderot -- Paris 7, France.\\ \printead{e2}}
\end{aug}

\received{\smonth{10} \syear{2013}}
\revised{\smonth{5} \syear{2014}}

\begin{abstract}
We present sharp tail asymptotics for the density and the
distribution function of linear combinations of correlated log-normal
random variables, that is, exponentials of components of a correlated
Gaussian vector. The asymptotic behavior turns out to depend on the
correlation between the components, and the explicit solution is found
by solving a tractable
quadratic optimization problem. These results can be used either to
approximate the probability of tail events directly, or to construct
variance reduction procedures to estimate these probabilities by Monte
Carlo methods. In particular, we propose an efficient importance
sampling estimator for the left tail of the distribution function of
the sum of log-normal variables. 
As a corollary of the tail asymptotics, we compute the asymptotics of
the conditional law of a Gaussian random vector given a linear
combination of exponentials of its components. In risk management
applications, this finding can be used for the systematic construction
of stress tests, which the financial institutions are required to
conduct by the regulators.
We also characterize the asymptotic behavior of the Value at Risk for
log-normal portfolios in the case where the confidence level tends to one.
\end{abstract}

\end{frontmatter}

\section{Introduction}\label{SI}

The multivariate log-normal distribution is a widely used stochastic
model in natural and social sciences, and linear combinations of
correlated log-normal random variables arise in many applications. For
example, in wireless communications, the distribution of the total
interference power coming from several sources is often described by a
sum of log-normal variables, and in financial risk management, a linear
combination of correlated log-normal variables may represent the value
of a portfolio of assets.
Since the distribution of such linear combinations is not known in
explicit form, a considerable effort has been devoted to developing
asymptotic approximations. In particular, Asmussen and Rojas-Nandayapa
\cite{AR-N} characterized the behavior of the right tail of the
distribution function of the sum of correlated log-normals.
Their results can also be deduced from the more recent studies of the
tail behavior of sums of dependent subexponential random variables \cite{FR,GT}.
On the other hand, Gao, Xu and Ye \cite{GXY} computed the asymptotics
of the left tail of the sum of \emph{two} correlated log-normal variables.
However, beyond these two cases, the tail behavior of linear
combinations of log-normal variables is not well understood so far.

In this paper, we present an explicit characterization of the tail
asymptotics of the density and the distribution function of arbitrary
linear combinations of correlated log-normal random variables. We find
new dependence patterns, very different from those which have been
established for the right tail of the log-normal sum and more generally
for the right tail of the sum of subexponential random variables.
The principle of a ``single big jump'' does not hold: the asymptotic
behavior is no longer determined by the single component with the
fattest tail but depends on the correlation between the components. 

Our paper contains two types of results. Firstly, we compute the tail
asymptotics of the distribution function and the density of a linear
combination of log-normal variables. These results can be used either
to estimate the probability of tail events directly, or to construct
efficient variance reduction procedures to estimate these probabilities
by Monte Carlo method. In particular, we propose an importance sampling
estimator for the left tail of the distribution function of the
log-normal sum, which is \emph{logarithmically efficient} in the sense
of Asmussen and Glynn \cite{AG}. In risk management applications, our
asymptotic formulas can be used to evaluate the probability of large
portfolio losses within the multidimensional Black--Scholes model.
Secondly, as a corollary of the tail asymptotics, we compute the
asymptotic law of a Gaussian vector conditional on a linear combination
of exponentials of its components. This finding can be used for the
systematic construction of stress tests, which the financial
institutions are required to conduct by the regulators.

In the present paper we focus on the multidimensional log-normal
distribution in view of its importance for applications. However, we
expect that our findings and the techniques we develop will stimulate
further studies of the multidimensional stochastic models and settings
in which the tail behavior is determined by the entire dependence
structure rather than by a single component.


\subsection*{Review of relevant literature}

The history and the applications of log-normal distributions are
reviewed in \cite{AB,CS,D,JK,LSA}. Sums and integrals of log-normal
variables and processes also play an important role in theoretical
probability and theoretical physics in relation to the Gaussian
multiplicative chaos \cite{K,RV}. A~considerable effort has been
devoted to the numerical approximations of the distribution function of
the sum of log-normal variables. In
\cite{BK,BKG}, the authors find approximations
to the density of the sum of log-normals based on approximations of the
characteristic function of the univariate log-normal distribution. A
similar path is taken in the papers \cite{ST} and \cite{TS} of
Senarante and Tellambura to develop deterministic numerical techniques
for the computation of the distribution function of the
log-normal sum.

A related approach is to bound the density of the log-normal sum from
above and below with some more or less easily computable expressions. Tellambura
\cite{T} provides bounds for the distribution function of a sum of 2 or
3 correlated log-normals, and also for the sum of any number of
equally-correlated log-normals. The paper \cite{VCDGHK} of Vanduffel
\textit{et al.} is devoted to approximations of the distribution function of the
log-normal sum by the distribution function of the conditional
expectation of such a sum with respect to an auxiliary conditioning
random variable. Another stream of literature discusses bounds for
tails of functions of general random vectors with fixed marginals (see
\cite{EP4} and references therein).

Motivated, in particular, by the applications in risk management,
several authors have studied the tail behavior of sums of log-normal
variables. As already mentioned, Asmussen and Rojas-Nandayapa \cite{AR-N} (see also the dissertation of Rojas-Nandayapa \cite{R-N})
characterized the behavior of the right tail of the distribution
function of the sum of log-normals. These results were used in Asmussen
\textit{et al.} \cite{ABJR-N}, to construct importance sampling Monte Carlo
estimators for the distribution function of the sum of log-normals in
the right tail. Understanding the left tail of the log-normal sum
turned out to be considerably more difficult. Szyszkowicz and
Yanikomeroglu \cite{SY} propose to approximate the left tail of the sum
of uncorrelated log-normal variables by a one-dimensional log-normal
distribution function. More recently, an important progress was made in
the article by Gao, Xu and Ye \cite{GXY} where explicit asymptotics of
the left tail of the distribution of the sum of two correlated
log-normals are presented. For a subclass of covariance matrices (see
Remark~\ref{gxy.rem} below) these authors also characterize the
asymptotic behavior of the left tail of the density of the sum of an
arbitrary number of log-normal variables.

The vast majority of publications discuss the sum of log-normal
variables. Linear combinations of such variables with coefficients of
different signs have received relatively little attention despite their
importance for applications, for instance, to spread option pricing in
finance (see Carmona and Durrleman \cite{CD}). One of the few exceptions is
the paper by Lo \cite{Lo} who considers the distribution of the sum and
the difference of two log-normal processes (geometric Brownian motions)
and presents a small-time approximation to these distributions.
Small-time and small-noise asymptotics of sums and differences of
geometric Brownian motions are also discussed in several papers dealing
with basket and spread option pricing \cite{ABBF,BL,BFL}.



The log-normal distribution is an example of a subexponential
distribution (see \cite{FKZ} for the definition
of subexponentiality). Numerous publications were devoted to tail
estimates for sums or more general functions of dependent
sub-exponential random variables (see \cite{AAK,FR,GT,KT,KA,YY} and the
references therein). The right-tail behavior of the difference of a
positive subexponential random variable and a dependent positive random
variable is
studied in \cite{AAK1}. 
Our paper focuses on the cases which cannot be dealt with using the
theory of subexponential distributions, such as the \emph{left} tail of
the sum of log-normals, or the right tail of the weighted sum with
weights of different signs.

\subsection*{Main notation}
Throughout the paper, we use boldface for
denoting vectors. In particular, $\mathbf{1}$ denotes the vector of
suitable dimension, all components of which are equal to $1$.
A strictly positive random vector $\ve{X}=(X_1,\ldots,X_n)$
such that the vector $\ve{Y}=(Y_1,\ldots,Y_n)$ with $Y_i=\log X_i$,
$1\le i\le n$, has an $n$-dimensional normal distribution
with mean vector $\veee{\mu}=(\mu_1,\ldots,\mu_n)$ and covariance
matrix~$\cm$, is called an $n$-dimensional log-normal vector with parameters
$\veee{\mu}$ and $\cm$. The elements of the matrix $\cm$ will be denoted by
$b_{ij}$. The distribution of $\ve{X}$ is called the $n$-dimensional
log-normal distribution and denoted by $\Lambda(\veee{\mu},\cm)$. In the
present paper, we make the standing assumption that $\llvert \cm\rrvert
>0$. The inverse matrix of the covariance matrix will be denoted by $\cm
^{-1}$, its elements will be denoted by $a_{ij}$, and we put $A_k=\sum_{j=1}^na_{kj}$, $1\le k\le n$.\vadjust{\goodbreak}
The log-normal distribution $\Lambda(\veee{\mu},\cm)$ admits a density defined by
%
\begin{eqnarray}
&& d^{\log}(x_1,\ldots,x_n)
\nonumber
\\[-8pt]
\label{Emult1}
\\[-8pt]
\nonumber
&& \quad =\frac{1}{(2\uppi)^{{n}/{2}}\sqrt{\llvert \cm\rrvert }x_1\cdots x_n} \exp \Biggl\{-\frac{1}{2}\sum
_{ i, j= 1}^n a_{ij} (\log x_i-
\mu_i) (\log x_j-\mu_j) \Biggr\},
\end{eqnarray}
where $x_i> 0$ for $1\le i\le n$. In particular, the one-dimensional
log-normal density with mean $\mu\in\mathbb{R}$ and variance $\sigma^2$
is given by
%
\begin{equation}\label{Emultf}
d^{\log}(x)=\frac{1}{\sqrt{2\uppi}\sigma x}\exp \biggl\{-\frac{1}{2\sigma^2} (\log x-
\mu)^2 \biggr\},\qquad  x> 0.
\end{equation}

For every integer $m$
with $1\le m\le n$, we consider the random variable
%
\begin{equation}\label{Ebass}
X^{(m)} = \sum_{k=1}^m
\mathrm{e}^{Y_k} - \sum_{k=m+1}^n
\mathrm{e}^{Y_k}.
\end{equation}
The support of $X^{(m)}$ is equal to $\mathbb{R}$ for $m=1,\ldots,n-1$,
and to $\mathbb{R}_+$ for $m=n$. For $m=n$, the variable $X^{(n)}$ is denoted
simply by $X$. The symbols $p^{(m)}$ and $p$ will stand for the density
of $X^{(m)}$
and $X$, respectively.

Our main goal in this paper is to characterize the
tail behavior of the distribution of the random variable $X^{(m)}$. We
are mainly interested in the asymptotic behavior of the right tail of
the variables $X^{(m)}$, $1\leq m\leq
n-1$, that is, the behavior of $\mathbb{P}[X^{(m)}>x]$ and $p^{(m)}(x)$
as $x\to\infty$, as well as the behavior of the left tail of $X$ (as
$x\to0$). The
right tail behavior of the distribution function of $X$ was completely
characterized in \cite{AR-N}, while in \cite{GXY}, a similar
characterization was obtained for the density.
The left tail behavior of $X^{(m)}$, $1\leq m\leq
n-1$, can be deduced from that of the right tail by exchanging the
signs of the
variables. Since positive coefficients can be incorporated into the mean
vector of ${Y}$, our results provide a complete characterization of the
tail behavior of a linear combination $\sum_{i=1}^n \lambda_i \mathrm{e}^{Y_i}$,
$\lambda_i\in\mathbb{R}$, of components of a log-normal
random vector $(\mathrm{e}^{Y_1},\ldots,\mathrm{e}^{Y_n})$.

\subsection*{Overview of the paper}
Section~\ref{Snachalo} deals with the left tail asymptotics of sums of
log-normal variables. This section is split into Section~\ref{SSfmr1}, where we formulate and discuss some of our main results, and
Section~\ref{SSproofs1}, which contains the proofs.

In Section~\ref{SSfmr1}, we formulate asymptotic formulas for the
distribution function and the distribution density in the general case,
under rather mild nondegeneracy conditions (Theorem~\ref{gen.thm} and
Corollary~\ref{Ccorr}). This is done by relating the tail asymptotics
to the quadratic optimization problem
%
\begin{equation}\label{qpintro.eq}
\min_{\ve{w}\in\Delta_n} \ve{w}^\perp\cm\ve{w},
\end{equation}
where $\Delta_n$ is the set of vectors in $\mathbb{R}^n$ whose
components are all non-negative and sum up to one. In particular, the
leading term in the asymptotics for both the distribution function and
the density is given by
%
\begin{equation}\label{leading.eq}
\exp \biggl\{-\frac{\log^2 x}{2 \min_{\ve{w}\in\Delta_n} \ve{w}^\perp\cm
\ve{w}} \biggr\}.
\end{equation}
This is in sharp contrast with the findings of \cite{AR-N} for the
asymptotics of the \emph{right tail} of the sum of log-normals, where
the leading term is
\[
\exp \biggl\{-\frac{\log^2 x}{2 \max_{i=1,\ldots,n} \cm_{ii}} \biggr\}.
\]
%
As an application of Theorem~\ref{gen.thm} and Corollary~\ref{Ccorr},
we characterize the asymptotic behavior of conditional Laplace
transforms of multidimensional Gaussian vectors (see Corollaries~\ref{laplace.cor} and~\ref{laplacedens.cor}). These estimates are used in
Section~\ref{rm.sec}, which deals with
stress testing of log-normal portfolios. Section~\ref{SSproofs1}
contains the proofs of the results formulated in
Section~\ref{SSfmr1}. We first establish asymptotic formulas for
the distribution function and the density in the special case (see
Lemma~\ref{easy.prop}), when the row sums of the inverse covariance
matrix are all strictly positive. Here we can apply Laplace's method to
the integral, characterizing the distribution density of the log-normal sum.
Lemma~\ref{easy.prop} is used in the proof of Theorem~\ref{gen.thm}.

In Section~\ref{Stail1}, the results obtained in Section~\ref{Snachalo} are extended to the case of the difference of two
log-normal sums (see the random variable in (\ref{Ebass})). These
extensions are not trivial, and they are new even in simple cases. We
find sharp asymptotic formulas for the right tails of distributions of
such differences. The proofs of the results obtained in Section~\ref{Stail1} are similar to those in Section~\ref{Snachalo}, but the
details are more complicated. In Section~\ref{Srta}, we formulate
several assertions concerning the tail asymptotics of log-normal differences.
Theorem~\ref{Tpresent} characterizes the asymptotic behavior of the
right tail of a log-normal difference.
In addition, the asymptotic behavior of the conditional Laplace
transform is characterized in Corollary~\ref{laplace2.cor}.
Section~\ref{SStrap2} is devoted to the proofs of these results.
Here we start with a special case, where Laplace's method can be
applied directly (see Lemma~\ref{Tstandard}), and reduce the general
case to the special one using quadratic programming methods.

In Section~\ref{num.sec}, we analyze the performance of our asymptotic
formulas via numerical examples, by comparing the theoretical results
with Monte Carlo computations. The convergence turns out to be quite
slow, which is consistent with logarithmic error bounds in our main
results. However, the asymptotic formulas provide a good order of
magnitude approximation for a wide range of values of $x$. This fact
enables us to design an importance sampling technique for evaluation of
the tail event probabilities by Monte Carlo method, which is
logarithmically efficient in the sense of Asmussen and Glynn \cite{AG}.

The last part of the paper (Section~\ref{rm.sec}) considers
applications of our asymptotic formulas to risk management in the
context of the multidimensional Black--Scholes model. This model, which
represents stock prices as exponentials of correlated Brownian motions,
remains widely used for the analysis of large portfolios. 
Our asymptotic theory provides two types of insights. First, it allows
to quantify the tail behavior of portfolios of log-normal stocks. For
example, for portfolios with positive weights, the leading term of the
probability of a large downside move is given by \eqref{leading.eq}.
This means that \eqref{qpintro.eq} measures the risk of the portfolio
in a downturn. Second, it provides better understanding of the behavior
of individual assets under various adverse scenarios. For instance, we
consider a typical stress scenario when the normalized value of a
benchmark portfolio (or index) drops to $x$ with $x$ small.
Theorem~\ref{stress.prop}
characterizes the asymptotic behavior of conditional expectations of
the individual assets in the original portfolio under such an adverse
scenario, when the benchmark portfolio has only positive weights. This
theorem shows that the assets in a market can be categorized into two
classes: those assets for which conditional expectations decay
proportionally to $x$ as $x\rightarrow0$ (safe assets), and those
assets, for which conditional expectations decay much faster than $x$
(dangerous assets). The safe assets are exactly those which have
strictly positive weights in the solution to the quadratic programming
problem \eqref{qpintro.eq}.
Results such as Theorem~\ref{stress.prop} may be employed for systematic construction
of stress tests, which banks and investment firms are required to
conduct by the regulatory bodies. {Finally, in Section~\ref{Siv},
we characterize the asymptotic behavior of the Value at Risk for a
log-normal portfolio as the confidence level tends to one (see Theorem~\ref{Tvarf}).}

\section{Asymptotic behavior of the left tail of a log-normal sum}\label
{Snachalo}
The present section studies the left tail asymptotics of the random variable
$X=\sum_{i=1}^n\mathrm{e}^{Y_i}$. 

\subsection{Left tail of a log-normal sum. Results and
discussions}\label{SSfmr1}
Denote by $\Delta_n$ the $n$-dimensional simplex defined by
\[
\Delta_n:= \Biggl\{\ve{w}\in\mathbb{R}^n\dvt
w_i\geq0,i=1,\ldots,n, \mbox{and }\sum_{i=1}^n
w_i = 1\Biggr\}
\]
and let $\bar{\ve{w}}\in\Delta_n$ be the unique vector such that
%
\begin{equation}\label{minprob}
\bar{\ve{w}}^\perp\cm\bar{\ve{w}} = \min_{\ve{w}\in\Delta_n}
\ve{w}^\perp\cm\ve{w}.
\end{equation}
The existence and uniqueness of $\bar{\ve{w}}$ follows from the
non-degeneracy of the matrix
$\cm$.
In the case where $A_k >0$ for $k=1,\ldots,n$,
%
\begin{equation}\label{easycase.eq}
\bar{\ve{w}} = \frac{\cm^{-1} \mathbf{1}}{\mathbf{1}^{\perp}
\cm^{-1}\mathbf{1}}\quad  \iff \quad  \bar{w}_k =
\frac{A_k}{\sum_{i=1}^n
A_i},\qquad k=1,\ldots,n,
\end{equation}
which means that $\bar{w}_k>0$ for $k=1,\ldots,n$.
In the general case, we
let
%
\begin{eqnarray}
\bar{n} &:=&  \operatorname{Card} \{i=1,\ldots,n\dvt  \bar{w}_i \neq0\},
\nonumber
\\[-8pt]
\label{nbar.eq}
\\[-8pt]
\nonumber
\bar I &:=&  \{i=1,\ldots,n \dvt  \bar{w}_i \neq0\}:= \bigl\{\bar{k}(1),\ldots,
\bar{k}(\bar{n})\bigr\},
\end{eqnarray}
$\bar{\veee{\mu}}\in\mathbb{R}^{\bar{n}}$ with $\bar{\mu}_i = \mu_{\bar{k}(i)}$,
and $\bar{\cm}\in M_{\bar{n}}(\mathbb{R})$ with $\bar\cm_{ij} = \cm_{\bar
k(i),\bar k(j)}$. The inverse matrix of $\bar\cm$ is denoted by $\bar
\cm^{-1}$ and its elements and row sums by $\bar a_{ij}$ and $\bar A_k
:= \sum_{j=1}^{\bar{n}}
\bar a_{kj}$.

In the present subsection, we are only dealing with the sum of the
exponentials of the random variables
$Y_1,\ldots,Y_n$. Since these variables are exchangeable, we can assume
with no loss of generality that for the covariance matrix $\cm$, $\bar
I = \{1,\ldots,\bar{n}\}$ with $\bar{n}\leq n$.
By the strict convexity of the objective function, the minimizer of
$ \min_{\ve{w}\in\Delta_{\bar{n}}} \ve{w}^\perp\bar\cm\ve{w}$ coincides with the first $\bar{n}$ components of $\bar{\ve{w}}$ and
therefore belongs to the
interior of the set $\mathbb{R}_+^{\bar{n}}$. The minimizer
over $\Delta_{\bar{n}}$ then coincides with the minimizer over the set
$\{\ve{w}\in
\mathbb{R}^{\bar{n}} \dvt \sum_{i=1}^{\bar{n}} w_i = 1\}$, which means that
\[
(\bar{w}_{i})_{i=1,\ldots,\bar{n}} = \frac{\bar\cm^{-1} \mathbf
1}{\mathbf{1}^\perp
\bar\cm^{-1}\mathbf{1}},
\]
or, equivalently,
%
\begin{equation}\label{awstar}
\bar{w}_k = \frac{\bar A_k}{\sum_{i=1}^{\bar{n}} \bar A_i},\qquad  k=1,\ldots ,\bar{n}.
\end{equation}
Since
$\sum_{i=1}^{\bar{n}} \bar A_i>0$ (the matrix $\bar\cm^{-1}$ is
positive definite),
this implies that $\bar A_k>0$ for $k=1,\ldots,\bar{n}$. Equation \eqref
{awstar} also leads to the following useful formula:
\[
\bar{\ve{w}}^\perp\cm\bar{\ve{w}} = \frac{1}{\mathbf{1}^\perp\bar\cm
^{-1} \mathbf{1}} =
\frac{1}{\sum_{i=1}^{\bar{n}} \bar A_i}.
\]

The following assumption will be used in the sequel:
\begin{enumerate}[$(\mathcal A)$]
\item[$(\mathcal A)$] For every $i\in\{1,\ldots,n\}\setminus\bar I$,
\[
\bigl(\ve{e}^i - \bar{\ve{w}}\bigr)^\perp\cm\bar{\ve{w}}
\neq0,
\]
where $\ve{e}^i\in\mathbb{R}^n$ satisfies $e^i_j = 1$ if $i=j$ and
$e^i_j=0$ otherwise.
\end{enumerate}

\begin{remark}
Assumption $(\mathcal A)$ is equivalent to the following:
\[
\bigl(\ve e^i - \bar{\ve{w}}\bigr)^\perp\cm\bar{\ve{w}} >
0,
\]
for every $i\in\{1,\ldots,n\}\setminus\bar I$. Indeed, the gradient of
the minimization functional $\frac{1}{2}
\ve{w}^\perp\cm\ve{w}$ at the point $\bar{\ve{w}}$ is given by $\cm\bar
{\ve{w}}$, and for
$\varepsilon>0$ small enough, $\bar{\ve{w}} + (\ve{e}^i-\bar{\ve{w}})\varepsilon$
clearly belongs to $\Delta_n$. Therefore $(\ve{e}^i - \bar{\ve{w}})^\perp\cm
\bar{\ve{w}} < 0$ would contradict the fact that $\bar{\ve{w}}$ is the
minimizer.

Assumption $(\mathcal A)$ is a natural nondegeneracy condition for our
problem. The following straightforward equality gives a
relation between the optimization problem in \eqref{minprob} and a similar
problem without the normalization constraint:
%
\begin{equation}\label{nonorm}
\inf_{\ve{w}\in\Delta_n, r\geq0} \frac{r^2}{2}{\ve{w}}^\perp\cm{
\ve{w}} - r = \inf_{\ve v\in\mathbb{R}^n\dvt v_i \geq0, i=1,\ldots,n} \frac{1}{2} \ve
v^\perp\cm\ve v - \mathbf{1}^\perp\ve v.
\end{equation}
A minimizer $\bar{\ve v}$ of the right-hand side can therefore be
constructed from
the minimizer $\bar{\ve{w}}$ of \eqref{minprob} as follows:
\[
\bar{\ve v} = \frac{\bar{\ve{w}}}{\bar{\ve{w}}^\perp\cm\bar{\ve{w}}}.
\]
Now, introducing the vector ${\veee\lambda}\in\mathbb{R}^n$ of Lagrange
multipliers
for the positivity constraints on the right-hand side of
\eqref{nonorm}, we get the Lagrangian
\[
\frac{1}{2} \ve v^\perp\cm\ve v - \mathbf{1}^\perp\ve v
- {\veee\lambda}^\perp\ve v.
\]
At the extremum therefore,
$\cm\bar{\ve v} = \mathbf{1}+{\veee\lambda}$,
or in other words,
\[
\frac{\cm\bar{\ve{w}}}{\bar{\ve{w}}^\perp\cm\bar{\ve{w}}} = \mathbf 1+{\veee\lambda}.
\]
Therefore, assumption $(\mathcal A)$ simply states that for the
constraints, which are
saturated, the Lagrange multipliers are not equal to zero (since the constraints
are inequalities, this is equivalent to the strict positivity for
the multipliers). This is generally true, except when the solution
of the unconstrained problem belongs to the boundary of the
domain defined by the constraints.
\end{remark}

The next assertion provides a sharp asymptotic formula
with an error estimate for the distribution function of the random
variable $X$, under assumption $(\mathcal A)$. A similar
formula for the distribution density of $X$ will be formulated below
(see Corollary~\ref{Ccorr}).

\begin{theorem}\label{gen.thm}
Suppose assumption $(\mathcal A)$ holds. Then, as $x\to0$,
%
\begin{eqnarray}
\mathbb{P}[X \leq x] &=& \frac{C}{\bar A_1+\cdots+\bar A_{\bar{n}}} \biggl(\log\frac{1}{x}
\biggr)^{-({1+\bar{n}})/{2}}x^{\sum_{k=1}^{\bar{n}}\bar A_k (\log({\bar A_1+\cdots+\bar A_n})/{\bar A_k}+\bar\mu_k )}
\nonumber
\label{Egenerass}\\[-8pt]
\\[-8pt]
\nonumber
&&{}\times \exp \biggl\{-\frac{1}{2}(\bar A_1+\cdots+\bar
A_{\bar{n}})\log ^2\frac{1}{x} \biggr\} \biggl(1+\mathrm{O}
\biggl( \biggl(\log\frac{1}{x} \biggr)^{-1} \biggr)
\biggr),
\end{eqnarray}
where
%
\begin{eqnarray}
C &=& \frac{1}{\sqrt{2\uppi}\sqrt{\llvert \bar\cm\rrvert }}\frac{\sqrt{\bar
A_1+\cdots+\bar A_{\bar{n}}}}{\sqrt{\bar A_1\cdots\bar A_{\bar{n}}}}
\nonumber
\\[-8pt]
\label{Efor1}
\\[-8pt]
\nonumber
&&{}\times\exp \Biggl\{-\frac{1}{2}\sum_{ i,j=1}^{\bar{n}}
\bar a_{ij} \biggl(\log\frac{\bar A_1+\cdots+\bar A_{\bar{n}}}{\bar
A_i}+\bar\mu_i
\biggr) \biggl(\log\frac{\bar A_1+\cdots+\bar A_{\bar{n}}}{\bar A_j}+\bar\mu _j \biggr) \Biggr\}.
\end{eqnarray}
\end{theorem}

\begin{remark}
Formula \eqref{Egenerass} can be rewritten in terms of the solution
$\bar{\ve{w}}$ to
the quadratic programming problem in
\eqref{minprob} as follows:
%
\begin{eqnarray}
\mathbb{P}[X \leq x] &=& \widetilde C \biggl(\log\frac{1}{x}
\biggr)^{-({1+\bar{n}})/{2}}\exp \biggl\{-\frac{(\log x-\bar{\ve{w}}^\perp\veee{\mu}
-\mathcal E(\bar{\ve{w}}))^2}{2 \bar{\ve{w}}^\perp\cm\bar{\ve{w}}} \biggr\}
\nonumber
\\[-8pt]
\label{genw.eq}
\\[-8pt]
\nonumber
&&{}\times \bigl(1+\mathrm{O} \bigl(|\log x|^{-1} \bigr) \bigr)
\end{eqnarray}
as $x\rightarrow0$, where $\mathcal E(\bar{\ve{w}}) = -\sum_{i=1}^{n}
\bar{w}_i \log\bar{w}_i$ and
\[
\widetilde C = C\bar{\ve{w}}^\perp\cm\bar{\ve{w}}\exp \biggl(
\frac
{\mathcal E(\bar{\ve{w}})^2}{2\bar{\ve{w}}^\perp\cm\bar{\ve{w}}} \biggr).
\]
The asymptotic behavior of the left tail of
the sum of log-normal random variables with positive coefficients is thus
intimately related to the quadratic programming problem formulated
in \eqref{minprob}. In particular, this problem determines which
components of the random vector influence the tail behavior.
\end{remark}

Theorem~\ref{gen.thm} and Corollary~\ref{Ccorr} below allow us to estimate
various conditional expectations. The next assertion provides a
characterization of
the limiting conditional law of the Laplace transform of $Y_1,\ldots,Y_n$, given that $X\leq x$.

\begin{corollary}\label{laplace.cor}
Suppose assumption $(\mathcal A)$ holds. Then, as $x\to0$, for any
$u\in\mathbb{R}^n$,
\begin{eqnarray*}
\mathbb E \bigl[\mathrm{e}^{\sum_{i=1}^n u_i Y_i}| X\leq x \bigr]& =& x^{\sum_{i=1}^n u_i \sum_{j=1}^{\bar{n}}
\bar A_j b_{ij}}
\\
&&{}\times\exp \Biggl\{ \sum_{i=1}^n
u_i \Biggl(\mu_i - \sum_{p,q=1}^{\bar{n}}
b_{pi} \bar a_{pq} \biggl(\log\frac{\bar A_1+\cdots+ \bar A_{\bar{n}}}{\bar A_q} + \bar
\mu_q \biggr) \Biggr) \Biggr\}
\\
&&{}\times\exp \Biggl\{\frac{1}{2}\sum_{i,j\notin\bar I}
u_i u_j \Biggl(b_{ij} - \sum
_{p,q=1}^{\bar{n}} \bar a_{pq}
b_{pi}b_{qj} \Biggr) \Biggr\} \biggl(1+\mathrm{O} \biggl( \biggl(
\log \frac{1}{x} \biggr)^{-1} \biggr) \biggr).
\end{eqnarray*}
\end{corollary}

\begin{remark}\label{remgt1}
Note that
\[
\sum_{j=1}^{\bar{n}} \bar A_j
b_{ij} = \frac{[\cm\bar
{\ve{w}}]_i}{ \bar{\ve{w}}^\perp\cm\bar{\ve{w}}} \cases{
=1, &\quad  $i \in\bar I$,\cr
>1, & \quad $i\notin\bar I$}
\]
by assumption $(\mathcal A)$. Let
%
\begin{equation}\label{Ela}
\bar\lambda_i = \frac{[\cm\bar
{\ve{w}}]_i}{ \bar{\ve{w}}^\perp\cm\bar{\ve{w}}} - 1.
\end{equation}
Then, Corollary~\ref{laplace.cor} implies that the conditional
distribution of the vector
\[
\ve Y - (\ve1+\bar{\veee\lambda})\log x,
\]
given $X\leq
x$, converges weakly to the (degenerate) Gaussian law with mean
\[
\mu^{\prime}_i = \mu_i - \sum
_{p,q=1}^{\bar{n}}b_{pi} \bar a_{pq}
\biggl(\log\frac{\bar A_1+\cdots+ \bar A_{\bar{n}}}{\bar A_q} + \bar\mu_q \biggr)
\]
and covariance matrix $\cm^{\prime}=(b^{\prime}_{ij})$ where
\[
b^{\prime}_{ij} = \Biggl\{b_{ij} - \sum
_{p,q=1}^{\bar{n}} \bar a_{pq}
b_{pi}b_{qj} \Biggr\}\mathbf 1_{i,j \notin\bar I}.
\]
Note that for $i\in\bar I$, the expression for $\mu^{\prime}_i$
simplifies to
\[
\mu^{\prime}_i = \log\frac{\bar A_i}{\bar A_1+\cdots+ \bar A_{\bar{n}}} = \log\bar{w}_{i}.
\]
\end{remark}

The next statement concerns the asymptotics of the distribution density
$p$ of the random variable $X$.

\begin{corollary}\label{Ccorr}
Suppose assumption $(\mathcal A)$ holds. Then, as $x\to0$,
%
\begin{eqnarray}
p(x)& =& C \biggl(\log\frac{1}{x} \biggr)^{({1-\bar{n}})/{2}} x^{-1+\sum_{k=1}^{\bar{n}}
\bar{A}_k (\log({\bar{A}_1+\cdots+\bar{A}_{\bar{n}}})/{\bar
A_k}+\bar{\mu}_k )}
\nonumber
\\[-8pt]
\label{Eforr1}
\\[-8pt]
\nonumber
&&{}\times\exp \biggl\{-\frac{1}{2}(\bar A_1+\cdots+\bar
A_{\bar{n}})\log ^2\frac{1}{x} \biggr\} \biggl(1+\mathrm{O}
\biggl( \biggl(\log \frac{1}{x} \biggr)^{-1} \biggr) \biggr),
\end{eqnarray}
where the constant $C$ is given by (\ref{Efor1}).
\end{corollary}

Corollary~\ref{Ccorr} implies that the conditional expectation in Corollary~\ref{laplace.cor} can be taken with respect to the event $\{X = x\}$.

\begin{corollary}\label{laplacedens.cor}
Suppose assumption $(\mathcal A)$ holds. Then, as $x\to0$, for any
$u\in\mathbb{R}^n$,
%
\begin{eqnarray}
\mathbb E\bigl[\mathrm{e}^{\sum_{i=1}^n u_iY_i}| X=x\bigr] &=& x^{\sum_{i=1}^nu_i\sum
_{j=1}^{\bar{n}} \bar A_jb_{ij}}
\nonumber
\\
\label{ELT} && {}\times \exp \Biggl( \sum_{i=1}^n
u_i \Biggl\{\mu_i - \sum_{p,q=1}^{\bar{n}}
b_{pi} \bar a_{pq} \biggl(\log\frac{\bar A_1+\cdots+ \bar A_{\bar{n}}}{\bar A_q} + \bar
\mu_q \biggr) \Biggr\} \Biggr)
\\
&&{}\times\exp \Biggl(\frac{1}{2}\sum_{i,j\notin\bar I}
u_i u_j \Biggl\{ b_{ij} - \sum
_{p,q=1}^{\bar{n}} \bar a_{pq} b_{pi}
b_{qj} \Biggr\} \Biggr) \biggl(1+\mathrm{O} \biggl( \biggl(\log
\frac{1}{x} \biggr)^{-1} \biggr) \biggr).\nonumber
\end{eqnarray}
\end{corollary}

\subsection*{Example: the sum of two log-normal variables}
Let $n = 2$, and denote the elements of the matrix $\cm$ by $b_{11} =
\sigma_1^2$, $b_{22} = \sigma_2^2$ and
$b_{12} = b_{21}=\rho\sigma_1\sigma_2$. To fix the ideas we assume
$\sigma_1\geq
\sigma_2$. Then, $\ve{w} = (v,1-v)^\perp$ and
\[
\ve{w}^\perp\cm\ve{w} = \sigma_1^2
v^2 + \sigma_2^2 (1-v)^2 + 2\rho
\sigma_1\sigma_2 v(1-v).
\]
%
Therefore, the solution to problem \eqref{minprob} is given by
\[
\bar{\ve{w}} =(\bar v, 1-\bar v)^\perp \qquad \mbox{with } \bar v=
\frac{\sigma_2(\sigma_2 -
\rho\sigma_1)}{\sigma_1^2 + \sigma_2^2 -2\rho\sigma_1 \sigma_2} \vee0
\]
and we have the following three cases:
\begin{itemize}
\item If $\rho< \frac{\sigma_2}{\sigma_1}$, then both weights are
strictly positive, assumption $(\mathcal A)$ holds, and the asymptotic
behavior of the density $p$ is as follows:
\begin{eqnarray*}
p(z) &=& \frac{C}{z\sqrt{|\log z|}}
\\
&&{}\times\exp\biggl\{-\frac{1}{2} \bigl(\mu_1 + x^* - \log z,
\mu_2 + y^* - \log z\bigr) \\
&&\quad\hspace*{20pt}{}\times\cm^{-1}\bigl(\mu_1 +
x^* - \log z, \mu_2 + y^* - \log z\bigr)^{\perp}\biggr\}
\\
&&{}\times \biggl(1+ \mathrm{O} \biggl(\frac{1}{|\log z|} \biggr) \biggr)
\end{eqnarray*}
as $z\rightarrow0$, with
\begin{eqnarray*}
C &= & \sqrt{\frac{\sigma_1^2 + \sigma_2^2 - 2\rho\sigma_1 \sigma_2}{2\uppi
(\sigma_2^2 - \rho\sigma_1 \sigma_2)(\sigma_1^2 - \rho\sigma_1
\sigma_2 )}},
\\
x^* &=&  \log\frac{\sigma_1^2 + \sigma_2^2 - 2\rho\sigma_1
\sigma_2}{\sigma_2^2 - \rho\sigma_1 \sigma_2} \quad \mbox{and}\quad  y^* = \log\frac{\sigma_1^2 + \sigma_2^2 - 2\rho\sigma
_1 \sigma_2}{\sigma_1^2 - \rho\sigma_1 \sigma_2}.
\end{eqnarray*}
\item If $\rho> \frac{\sigma_2}{\sigma_1}$, then $\bar{\ve{w}} =
(0,1)^\perp$, assumption $(\mathcal A)$ holds, and the asymptotic
behavior of the density is characterized by
\[
p(z) = \frac{1}{z \sigma_2 \sqrt{2\uppi}} \mathrm{e}^{ - {(\log z - \mu_2
)^2}/({2\sigma_2^2 })} \biggl(1+\mathrm{O} \biggl(\frac{1}{|\log z|}
\biggr) \biggr)
\]
as $z\rightarrow0$. Note that in this case the asymptotic behavior of
$p$ is determined by the second component only.
\item The case, where $\rho= \frac{\sigma_2}{\sigma_1}$, is
exceptional. Here we have $\bar{\ve{w}} = (0,1)^\perp$, but assumption
$(\mathcal A)$ does not hold. Thus, Theorem~\ref{gen.thm} can not be applied.
\end{itemize}
In \cite{GXY}, Gao, Xu, and Ye characterize the left tail behavior of
the sum of two log-normal variables in all the three cases described
above. It follows from the results established in \cite{GXY} that the
asymptotic behavior of the density $p$ in the exceptional case is
qualitatively different from the behavior of $p$ in the cases where
$\rho>\frac{\sigma_2}{\sigma_1}$ or $\rho< \frac{\sigma_2}{\sigma_1}$.
This shows that one can not relax assumption~$(\mathcal A)$ without
changing the form of the asymptotics, which means that, in a sense,
assumption~$(\mathcal A)$ is optimal.

\begin{remark}
In the second case of the above example, the variance of the second
component is so small, that it completely dominates the asymptotic
behavior of the left tail of the log-normal sum, so that, in a way, we
recover the law of one jump. When this is the case, the asymptotic
behavior of the distribution function can be characterized using an
elementary argument given in the following proposition. In the text of
the proposition and its proof, we denote $\sigma_i = \sqrt{\cm_{ii}}$
and $\rho_{ij} = \frac{\cm_{ij}}{\sigma_i\sigma_j}$ for $1\leq i,j\leq
n$. We thank the anonymous referee for bringing this argument to our attention.

\begin{proposition}
Assume that for some $i$, $\sigma_i < \rho_{ij} \sigma_j$, $\forall
j\neq i$. Then, as $x\to0$,
%
\begin{equation}\label{elemarg.eq}
\mathbb{P}\mathbb[X\leq x] = \frac{\sigma_i }{\log({1}/{x}) \sqrt
{2\uppi}} \mathrm{e}^{-{(\log x-\mu_i)^2}/({2\sigma_i^2})} \biggl(1+ \mathrm{O}
\biggl(\frac
{1}{\sqrt{\log{1}/{x}}} \biggr) \biggr).
\end{equation}
\end{proposition}

\begin{pf}
By Jensen's inequality,
\[
\mathbb E\bigl[(x-X)^+\bigr] \geq\mathbb E\bigl[\bigl(x-\mathbb
E[X|Y_i]\bigr)^+\bigr].
\]
An easy computation shows that
\[
\mathbb E[X|Y_i] = \mathrm{e}^{Y_i} + \sum
_{j\neq i} \exp \biggl(\mu_j + \sigma_j
\rho_{ij} (Y_i-\mu_i)/\sigma_i +
\frac{1}{2}\sigma_j^2 \bigl(1-
\rho_{ij}^2\bigr) \biggr) \leq \mathrm{e}^{Y_i} + c
\mathrm{e}^{\alpha Y_i}
\]
for some constants $c>0$ and $\alpha> 1$.
Combining this with a simple monotonicity argument, we get
\[
\mathbb E\bigl[\bigl(x-\mathrm{e}^{Y_i} - c \mathrm{e}^{\alpha Y_i}\bigr)^+\bigr]\leq
\mathbb E\bigl[(x-X)^+\bigr] \leq \mathbb E\bigl[\bigl(x-\mathrm{e}^{Y_i}\bigr)^+
\bigr].
\]
Since, on the event $\{x> \mathrm{e}^{Y_i} + c \mathrm{e}^{\alpha Y_i}\}$, $Y_i < \log
x$, we also have,
%
\begin{equation}\label{2side}
\mathbb E\bigl[\bigl(x-c x^\alpha- \mathrm{e}^{Y_i} \bigr)^+\bigr]\leq
\mathbb E\bigl[(x-X)^+\bigr] \leq \mathbb E\bigl[\bigl(x-\mathrm{e}^{Y_i}\bigr)^+
\bigr].
\end{equation}
By standard arguments,
\begin{eqnarray*}
\mathbb E\bigl[\bigl(x-\mathrm{e}^{Y_i}\bigr)^+\bigr] &=&  -\mathrm{e}^{\mu_i + \sigma_i^2/2} N
\biggl(-\frac{\mu
_i + \sigma_i^2 - \log x }{\sigma_i} \biggr) + x N \biggl(-\frac{\mu_i -
\log x }{\sigma_i} \biggr)
\\
& =& x
\sigma_i n \biggl(-\frac{\mu_i - \log x }{\sigma_i} \biggr) \biggl(
\frac{\sigma_i^2}{\log^2 {1}/{x}} + \mathrm{O} \biggl(\frac{1}{\log^3 {1}/{x}} \biggr) \biggr),\quad \mbox{as } x\to
0,
\end{eqnarray*}
where $N$ is the standard normal distribution function and $n$ is the
standard normal density. From the previous estimate and \eqref{2side},
we deduce that
\[
\mathbb E\bigl[(x-X)^+\bigr] = \frac{x \sigma_i^3}{\log^2 ({1}/{x}) \sqrt{2\uppi
}} \mathrm{e}^{-{(\mu_i-\log x)^2}/({2\sigma_i^2})} \biggl(1+ \mathrm{O}
\biggl(\frac{1}{\log{1}/{x}} \biggr) \biggr)\quad  \mbox{as }x\to 0.
\]
To obtain the asymptotics for the distribution function, we use the
following simple bound: for $\delta\in(0,x)$,
\[
\frac{\mathbb E[(x-X)^+]-\mathbb E[(x-\delta-X)^+]}{\delta}\leq\mathbb{P}[X\leq x ] \leq\frac{\mathbb E[(x+\delta-X)^+]-\mathbb
E[(x-X)^+]}{\delta}.
\]
Taking $\delta= \frac{x}{ (\log{1}/{x} )^{{3}/{2}}}$,
we then obtain formula \eqref{elemarg.eq} after some straightforwad
computations.
\end{pf}
\end{remark}

\subsection{Left tail of a log-normal sum. Proofs}\label{SSproofs1}

The following preliminary result characterizes the asymptotic behavior
of the distribution
function and the density of the random variable $X$ in the tail regime
in the special case where
$A_k>0$ for all $k=1,\ldots,n$, or, equivalently, $\bar{n} =n$.

\begin{lemma}\label{easy.prop}
Assume that $A_k>0, k=1,\ldots,n$. Then, as $x\rightarrow0$,
%
\begin{eqnarray}
p(x)&=& C \biggl(\log\frac{1}{x} \biggr)^{({1-n})/{2}}x^{-1+\sum
_{k=1}^nA_k (\log({A_1+\cdots+A_n})/{A_k}+\mu_k )}
\nonumber
\\[-8pt]
 \label{EPro3}
 \\[-8pt]
\nonumber
&&{}\times\exp \biggl\{-\frac{1}{2}(A_1+\cdots+A_n)
\log^2\frac{1}{x} \biggr\} \biggl(1+\mathrm{O} \biggl( \biggl(\log
\frac{1}{x} \biggr)^{-1} \biggr) \biggr),
\end{eqnarray}
and
%
\begin{eqnarray}
\mathbb{P}[X \leq x]&=& \frac{C}{A_1+\cdots+A_n} \biggl(\log\frac{1}{x}
\biggr)^{-({1+n})/{2}}x^{\sum_{k=1}^nA_k (\log({A_1+\cdots
+A_n})/{A_k}+\mu_k )}
\nonumber
\\[-8pt]
\label{EProf}\\[-8pt]
\nonumber
&&{}\times \exp \biggl\{-\frac{1}{2}(A_1+\cdots+A_n)
\log^2\frac{1}{x} \biggr\} \biggl(1+\mathrm{O} \biggl( \biggl(\log
\frac{1}{x} \biggr)^{-1} \biggr) \biggr),
\end{eqnarray}
where
%
\begin{eqnarray}
C &=& \frac{1}{\sqrt{2\uppi}\sqrt{\llvert \cm\rrvert }}\frac{\sqrt{A_1+\cdots
+A_n}}{\sqrt{A_1\cdots A_n}}
\nonumber
\\[-8pt]
\label{EPro4}\\[-8pt]
\nonumber
&&{}\times \exp \Biggl\{-\frac{1}{2}\sum_{ i,j=1}^n
a_{ij} \biggl(\log\frac
{A_1+\cdots+A_n}{A_i}+\mu_i \biggr) \biggl(
\log\frac{A_1+\cdots+A_n}{A_j}+\mu_j \biggr) \Biggr\}.
\end{eqnarray}
\end{lemma}

\begin{remark}\label{gxy.rem}
Formula \eqref{EPro3} for the density of the sum of log-normal random
variables under the assumption $A_k>0$ for $k=1,\ldots,n$ was given in
the paper \cite{GXY}, but with a very different notation and without
error estimates.
\end{remark}
\begin{pf*}{Proof of Lemma~\ref{easy.prop}}
We will first prove the formula in (\ref{EPro3}).
The distribution function of $X$ is given by
\[
\mathbb{P}[X<x] = \int_0^{x}\mathrm{d}y_1\int
_0^{x-y_1}\mathrm{d}y_2\cdots\int
_0^{x-y_1-y_2-\cdots-y_{n-1}} \mathrm{d}^{\log}(y_1,
\ldots,y_n)\,\mathrm{d}y_n.
\]
Differentiating the previous formula with respect to $x$ and making the
change of variables
%
\begin{eqnarray}
x_1 &=& y_1/x,\qquad  x_2 = y_2/x,\qquad
\ldots,\qquad
\nonumber
\\[-8pt]
\label{Edepend}
\\[-8pt]
\nonumber
 x_{n-1} &=& y_{n-1}/x,\qquad  x_n=1-(y_1+y_2+
\cdots+y_{n-1})/x,
\end{eqnarray}
we see that the density $p$ of the random variable $X$ can be
represented as follows:
%
\begin{equation}\label{Evis1}
p(x)=x^{n-1}\int_0^1\mathrm{d}x_1
\int_0^{1-x_1}\mathrm{d}x_2\cdots\int
_0^{1-x_1-\cdots-x_{n-2}} \mathrm{d}^{\log}(x x_1,
\ldots,x x_{n})\,\mathrm{d}x_{n-1}.
\end{equation}
Remark that $x_n$ is not an independent variable but a function of
$x_1,\ldots,x_{n-1}$.
Now, taking into account (\ref{Emult1}) and (\ref{Evis1}), we see
that for every $x> 0$,
%
\begin{eqnarray}
p(x)&=&\frac{1}{(2\uppi)^{{n}/{2}}\sqrt{\llvert \cm\rrvert }x} \int_0^1\mathrm{d}x_1
\int_0^{1-x_1}\mathrm{d}x_2\cdots \nonumber
\\[-8pt]
\label{maindens}
\\[-8pt]
\nonumber
&&{}\times\int
_0^{1-x_1-\cdots-x_{n-2}}\frac
{1}{x_1\cdots x_n}  \exp \Biggl\{-\frac{1}{2}\sum_{i,j=1}^n
a_{ij}\bigl(\log(xx_i)-\mu _i\bigr) \bigl(
\log(xx_j)-\mu_j\bigr) \Biggr\}\,\mathrm{d}x_{n-1}.
\end{eqnarray}

In the tail regime ($x\to0$), we can isolate the effect of $x$ in
formula \eqref{maindens}:
%
\begin{eqnarray}
p(x)&=& \frac{1}{(2\uppi)^{{n}/{2}}\sqrt{\llvert \cm\rrvert }x} \exp \biggl\{-\frac{1}{2}(A_1+
\cdots+A_n)\log^2\frac{1}{x} \biggr\}\int
_0^1\mathrm{d}x_1\int_0^{1-x_1}\mathrm{d}x_2
\cdots
\nonumber
\\[-8pt]
\label{Evis2}\\[-8pt]
\nonumber
&&{}\times\int_0^{1-x_1-\cdots-x_{n-2}} \Phi(x_1,
\ldots,x_{n-1})\exp \biggl\{-\log\frac{1}{x}\Psi(x_1,
\ldots ,x_{n-1}) \biggr\}\,\mathrm{d}x_{n-1},
\end{eqnarray}
where
%
\begin{equation}
\label{Evis3}\Phi(x_1,\ldots,x_{n-1})=\frac{1}{x_1\cdots x_n} \exp \Biggl\{-
\frac
{1}{2}\sum_{i,j=1}^n
a_{ij}\biggl(\log\frac{1}{x_i}+\mu_i\biggr) \biggl(
\log\frac
{1}{x_j}+\mu_j\biggr) \Biggr\}
\end{equation}
and
%
\begin{equation}\label{Evis4}
\Psi(x_1,\ldots,x_{n-1})=\sum_{k=1}^nA_k
\biggl(\log\frac{1}{x_k}+\mu _k \biggr).
\end{equation}

It is clear from formula (\ref{Evis2}) that it suffices to
characterize the asymptotic behavior as $\theta\to\infty$ of the integral
%
\begin{eqnarray}
I(\theta) &=& \int_0^1\mathrm{d}x_1\int
_0^{1-x_1}\mathrm{d}x_2\cdots
\nonumber
\\[-8pt]
\label{Evis5}\\[-8pt]
\nonumber
&&{}\times\int_0^{1-x_1-\cdots-x_{n-2}} \Phi(x_1,
\ldots,x_{n-1})\exp \bigl\{-\theta\Psi(x_1,
\ldots,x_{n-1}) \bigr\}\,\mathrm{d}x_{n-1}.
\end{eqnarray}
We will use the higher-dimensional extension of Laplace's method in the proof.

Recall that $A_k> 0$ for all $1\le k\le n$. We have
\[
\frac{\partial\Psi}{\partial x_k}=-\frac{A_k}{x_k} +\frac{A_n}{1-x_1-\cdots-x_k-\cdots-x_{n-1}}
\]
for all $1\le k\le n-1$. It follows that the function $\Psi$ has a
unique critical point
${x}^{*}=(x_1^{*},\ldots,x_{n-1}^{*})$ where
\[
x^{*}_k=\frac{A_k}{A_1+A_2+\cdots+A_n}
\]
for $1\le k\le n-1$.
Note that the critical point ${x}^{*}$ belongs to the interior of the
integration
set in (\ref{Evis5}), and moreover, this point is the global minimum
point of the function $\Psi$. Next, using formula (8.3.50) in \cite{BH}, we obtain
%
\begin{eqnarray}
I \biggl(\log\frac{1}{x} \biggr)&=& \frac{1}{\sqrt{\det(H({x}^{*}))}} \biggl(
\frac{2\uppi}{\log{1}/{x}} \biggr)^{({n-1})/{2}} \Phi\bigl(x_1^{*},
\ldots,x_{n-1}^{*}\bigr)
\nonumber
\\[-8pt]
\label{Evis6}
\\[-8pt]
\nonumber
&&{}\times \exp \biggl\{-\log\frac{1}{x}\Psi\bigl(x_1^{*},
\ldots,x_{n-1}^{*}\bigr) \biggr\} \biggl(1+\mathrm{O} \biggl( \biggl(
\log\frac{1}{x} \biggr)^{-1} \biggr) \biggr)
\end{eqnarray}
as $x\rightarrow0$, where $H({x}^{*})$ is the Hessian matrix of the
function $\Psi$ evaluated at the critical point
${x}^{*}$. Note that
%
\begin{eqnarray}
&& \Phi\bigl(x_1^{*},\ldots,x_{n-1}^{*}
\bigr)\nonumber\\
&& \quad = \frac{(A_1+\cdots+A_n)^n}{A_1\cdots
A_n}
\nonumber
\\[-8pt]
\label{Evis7}\\[-8pt]
\nonumber
&&\qquad {}\times\exp \Biggl\{-\frac{1}{2}\sum_{k=1}^na_{kk}
\biggl(\log\frac{A_1+\cdots
+A_n}{A_k}+\mu_k \biggr)^2
\\
&&\qquad \hspace*{31pt}{}-\sum_{1\le i< j\le n}a_{ij} \biggl(\log
\frac{A_1+\cdots
+A_n}{A_i}+\mu_i \biggr) \biggl(\log\frac{A_1+\cdots+A_n}{A_j}+
\mu_j \biggr) \Biggr\}\nonumber
\end{eqnarray}
and
%
\begin{equation}\label{Evis8}
\Psi\bigl(x_1^{*},\ldots,x_{n-1}^{*}
\bigr)=\sum_{k=1}^nA_k
\biggl(\log\frac
{A_1+\cdots+A_n}{A_k}+\mu_k \biggr).
\end{equation}
Moreover, since
\[
\frac{\partial^2\Psi}{\partial x_k^2}\bigl({x}^{*}\bigr)=(A_1+
\cdots+A_n)^2 \biggl(\frac{1}{A_k}+\frac{1}{A_n}
\biggr), \qquad 1\le k\le n-1,
\]
and
\[
\frac{\partial^2\Psi}{\partial x_i\,\partial x_j}\bigl({x}^{*}\bigr)=(A_1+\cdots
+A_n)^2\frac{1}{A_n},\qquad  1\le i< j\le n-1,\
\]
we have
%
\begin{eqnarray}
\det\bigl(H\bigl({x}^{*}\bigr)\bigr)& =&(A_1+
\cdots+A_n)^{2n-2}
\nonumber
\\[-8pt]
\label{Evis9}
\\[-8pt]
\nonumber
&&{}\times\left\vert
\begin{array}{cccc} A_1^{-1}+A_n^{-1}
& A_n^{-1} & \cdots& A_n^{-1}
\\
A_n^{-1} & A_2^{-1}+A_n^{-1}
& \cdots& A_n^{-1}
\\
\cdots& \cdots& \cdots& \cdots
\\
A_n^{-1} & A_n^{-1} &\cdots&
A_{n-1}^{-1}+A_n^{-1} \end{array}
\right\vert.
\end{eqnarray}
Next, using (\ref{Evis9}) and making long and tedious computations, we
get the following equality:
%
\begin{equation}\label{Evis10}
\det\bigl(H\bigl({x}^{*}\bigr)\bigr)=\frac{(A_1+A_2+\cdots+A_n)^{2n-1}}{A_1A_2\cdots A_n}.
\end{equation}

Finally, taking into account (\ref{Evis2}), (\ref{Evis6}), and (\ref{Evis7}),
(\ref{Evis8}), and (\ref{Evis10}), we complete the proof of formula
(\ref{EPro3}) in
Lemma~\ref{easy.prop}. Formula (\ref{EProf}) can be derived by
integrating formula
(\ref{EPro3}), or we can prove (\ref{EProf}) directly by employing
the same methods as those used in the proof of
(\ref{EPro3}).
\end{pf*}

\begin{pf*}{Proof of Theorem~\ref{gen.thm}}
Let $k\in\{\bar{n}, \ldots,n-1\}$, $z\in (0,\frac{1}{2} )$,
and $a,b$
be such that $\mathrm{e}^a + \mathrm{e}^b = z$. Then,
\begin{eqnarray*}
\mathbb{P}\bigl[\mathrm{e}^{Y_1}+\cdots+ \mathrm{e}^{Y_{k+1}} \leq z\bigr] &\geq &
\mathbb{P}\bigl[\mathrm{e}^{Y_1}+\cdots + \mathrm{e}^{Y_{k}} \leq \mathrm{e}^a,
Y_{k+1} \leq b\bigr]
\\
& = & \mathbb{P}\bigl[\mathrm{e}^{Y_1}+\cdots + \mathrm{e}^{Y_{k}} \leq
\mathrm{e}^a\bigr] - \mathbb{P}\bigl[\mathrm{e}^{Y_1}+\cdots + \mathrm{e}^{Y_{k}}
\leq \mathrm{e}^a, Y_{k+1} > b\bigr].
\end{eqnarray*}
Note that $k\geq\bar{n}$ implies that $\sum_{i=1}^{k}\bar{w}_i = 1$.
The second term in the above formula can be estimated as
follows:
\begin{eqnarray*}
\mathbb{P}\bigl[\mathrm{e}^{Y_1}+\cdots+ \mathrm{e}^{Y_{k}} \leq \mathrm{e}^a,
Y_{k+1} > b\bigr] &\leq & \mathbb{P}[\bar{w}_1 Y_1 +
\cdots+ \bar{w}_{k} Y_{k}  \leq  a, Y_{k+1} > b]
\\
& \leq & \mathbb{P}\bigl[\bar{w}_1 Y_1 + \cdots+ \bar{w}_{k} Y_{k} - a \leq \alpha( Y_{k+1} - b)\bigr]
\\
& =&  N \biggl(\frac{a-\alpha b - \mathbb E[\bar{\ve{w}}^{\perp} \ve Y -
\alpha
Y_{k+1}]}{\sqrt{\operatorname{Var}  [\bar{\ve{w}}^{\perp} \ve Y - \alpha
Y_{k+1}]}} \biggr)
\end{eqnarray*}
for every $\alpha>0$. Now, let
\[
x_{k+1} = \frac{(\ve{e}^{k+1})^\perp
\cm\bar{\ve{w}}}{\bar{\ve{w}}^\perp\cm\bar{\ve{w}}}>1
\]
(the inequality
follows from assumption $(\mathcal A)$) and choose
\[
b = \bigl(1 + (x_{k+1}-1)/2\bigr) \log z \quad \Rightarrow\quad  a = \log z + \log
\bigl(1-z^{({x_{k+1}-1})/{2}}\bigr).
\]
Noting that
\[
\operatorname{Var} \bigl[\bar{\ve{w}}^{\perp} \ve Y - \alpha Y_{k+1}
\bigr] = \bar{\ve{w}}^{\perp} \cm\bar{\ve{w}}(1 - 2 \alpha
x_{k+1}) + \alpha^2 \cm_{k+1,k+1},
\]
and making the above substitutions, we obtain, for $\alpha$ small enough,
\begin{eqnarray*}
&& \frac{a-\alpha b - \mathbb E[\bar{\ve{w}}^{\perp} \ve Y - \alpha
Y_{k+1}]}{\sqrt{\operatorname{Var}  [\bar{\ve{w}}^{\perp} \ve Y - \alpha
Y_{k+1}]}}
\\
&& \quad = \frac{ (1-\alpha({x_{k+1}+1})/{2}) \log z+ \log(1-z^{({x_{k+1}-1})/{2}}) - \bar
{\ve{w}}^\perp\veee{\mu}+ \alpha\mu_{k+1}}{\sqrt{\bar{\ve{w}}^{\perp}
\cm\bar{\ve{w}}}\sqrt{1-2\alpha x_{k+1} + {\alpha^2
\cm_{k+1,k+1}}/({\bar{\ve{w}}^{\perp} \cm\bar{\ve{w}}})}}
\\
&& \quad \leq \frac{\log z}{\sqrt{\bar{\ve{w}}^{\perp} \cm\bar{\ve{w}}}} \frac{1-\alpha
({x_{k+1}+1})/{2}}{\sqrt{1-2\alpha x_{k+1} + {\alpha^2\cm_{k+1,k+1}}/({\bar{\ve{w}}^{\perp} \cm\bar{\ve{w}}})}}+ C_{k+1},
\end{eqnarray*}
where $C_{k+1}$ is a constant which does not depend on $z$. Next, for
$\alpha$
small enough,
\begin{eqnarray*}
\frac{1-\alpha
({x_{k+1}+1})/{2}}{\sqrt{1-2\alpha x_{k+1} + {\alpha^2\cm_{k+1,k+1}}/({\bar{\ve{w}}^{\perp} \cm\bar{\ve{w}})}}} & \geq  &
\biggl(1-\alpha \frac{x_{k+1}+1}{2} \biggr) \sqrt{1+2
\alpha x_{k+1} - \frac{\alpha^2
\cm_{k+1,k+1}}{\bar{\ve{w}}^{\perp} \cm\bar{\ve{w}}}}
\\
& \geq & \biggl(1-\alpha \frac{x_{k+1}+1}{2} \biggr) \biggl(1+\alpha
x_{k+1} - \frac{\alpha^2
\cm_{k+1,k+1}}{2\bar{\ve{w}}^{\perp} \cm\bar{\ve{w}}} \biggr)
\\
& =&  1 + \alpha\frac{x_{k+1}-1}{2} + \mathrm{O}\bigl(\alpha^2\bigr), \qquad \alpha
\to0.
\end{eqnarray*}
Now it is easy
to see that by choosing $\alpha$ small enough, one can always find
$\varepsilon_{k+1} >0$ such that
\[
\frac{a-\alpha b - \mathbb E[\bar{\ve{w}}^{\perp} \ve Y - \alpha
Y_{k+1}]}{\sqrt{\operatorname{Var}  [\bar{\ve{w}}^{\perp} \ve Y - \alpha
Y_{k+1}]}} \leq\frac{\log z}{\sqrt{\bar{\ve{w}}^{\perp} \cm\bar{\ve{w}}}} (1+ \varepsilon_{k+1}) +
C_{k+1}.
\]
We conclude that
\begin{eqnarray*}
\mathbb{P}\bigl[\mathrm{e}^{Y_1}+\cdots+\mathrm{e}^{Y_{k+1}}\leq z\bigr]& \geq &
\mathbb{P}\bigl[\mathrm{e}^{Y_1}+\cdots+\mathrm{e}^{Y_{k}}\leq z\bigl(1-z^{({x_{k+1} - 1})/{2}}
\bigr)\bigr]
\\
&&{}- N \biggl(\frac{\log z}{\sqrt{\bar{\ve{w}}^\perp\cm\bar
{\ve{w}}}}(1+\varepsilon_{k+1}) + C_{k+1}
\biggr).
\end{eqnarray*}

Let us first apply this formula for $k=\bar{n}$ and $z=x$. Since
\[
N(y) = \mathrm{O} \biggl(\frac{\mathrm{e}^{-{y^2}/{2}}}{|y|} \biggr)
\]
as $y\to-\infty$, and
\[
\bar{\ve{w}}^\perp\cm\bar{\ve{w}} = \frac{1}{\sum_{i=1}^{\bar{n}} \bar A_i},
\]
using Lemma~\ref{easy.prop} to compute the asymptotics of
\[
\mathbb{P}\bigl[\mathrm{e}^{Y_1}+\cdots+\mathrm{e}^{Y_{\bar{n}}}\leq x\bigr]
\]
we have that
\[
\frac{N ({\log x}/{\sqrt{\bar{\ve{w}}^\perp\cm\bar
{\ve{w}}}}(1+\varepsilon_{\bar{n}+1}) + C_{\bar{n}+1} ) }{\mathbb{P}[\mathrm{e}^{Y_1}+\cdots+\mathrm{e}^{Y_{\bar{n}}}\leq x] } = \mathrm{O} \biggl( \biggl(\log\frac
{1}{x}
\biggr)^{-1} \biggr)
\]
as $x\to0$. On the other hand, by Lemma~\ref{easy.prop}, for $\delta>0$
\begin{eqnarray*}
&&\mathbb{P}\bigl[\mathrm{e}^{Y_1}+\cdots+\mathrm{e}^{Y_{\bar{n}}}\leq x
\bigl(1-x^{\delta}\bigr)\bigr]
\\
&&\quad = \frac{C}{\bar
A_1+\cdots+\bar A_{\bar{n}}} \biggl(\log\frac{1}{x} -\log\bigl(1-x^\delta
\bigr) \biggr)^{-{\bar{n}}/{2}} \bigl(x\bigl(1-x^\delta\bigr)
\bigr)^{\sum_{k=1}^{\bar{n}} \bar A_k  (\log({\bar
A_1+\cdots+\bar A_{\bar{n}}})/{\bar A_k}+\bar
\mu_k )}
\\
&&\qquad{}\times\exp \biggl(-\frac{1}{2}(\bar A_1+\cdots+\bar A_{\bar{n}}
) \biggl(\log\frac{1}{x} -\log\bigl(1-x^\delta\bigr)
\biggr)^2 \biggr) \biggl(1+\mathrm{O} \biggl( \biggl(\log \frac{1}{x}
\biggr)^{-1} \biggr) \biggr)
\\
&&\quad  = \frac{C}{\bar
A_1+\cdots+\bar A_{\bar{n}}} \biggl(\log\frac{1}{x} \biggr)^{-{\bar{n}}/{2}}
x^{\sum_{k=1}^{\bar{n}} \bar A_k  (\log
({\bar{A}_1+\cdots+\bar A_{\bar{n}}})/{\bar A_k}+\bar
\mu_k )}
\\
&&\qquad {}\times\exp \biggl(-\frac{1}{2}(\bar A_1+\cdots+\bar A_{\bar{n}}
) \log^2\frac{1}{x} \biggr) \biggl(1+\mathrm{O} \biggl( \biggl(\log
\frac{1}{x} \biggr)^{-1} \biggr) \biggr)
\end{eqnarray*}
as $x\rightarrow0$. Therefore, we have shown that
\[
\mathbb{P}\bigl[\mathrm{e}^{Y_1} + \cdots+ \mathrm{e}^{Y_{\bar{n}+1}} \leq z\bigr] \geq
\mathbb{P}\bigl[\mathrm{e}^{Y_1} + \cdots+ \mathrm{e}^{Y_{\bar{n}}} \leq z\bigr] \biggl(1+\mathrm{O}
\biggl( \biggl(\log \frac{1}{x} \biggr)^{-1} \biggr) \biggr)
\]
as $x\rightarrow0$, and since clearly,
\[
\mathbb{P}\bigl[\mathrm{e}^{Y_1} + \cdots+ \mathrm{e}^{Y_{\bar{n}+1}} \leq z\bigr] \leq
\mathbb{P}\bigl[\mathrm{e}^{Y_1} + \cdots+ \mathrm{e}^{Y_{\bar{n}}} \leq z\bigr],
\]
we also get
\[
\mathbb{P}\bigl[\mathrm{e}^{Y_1} + \cdots+ \mathrm{e}^{Y_{\bar{n}+1}} \leq z\bigr] = \mathbb{P}\bigl[\mathrm{e}^{Y_1} + \cdots+ \mathrm{e}^{Y_{\bar{n}}} \leq z\bigr] \biggl(1+\mathrm{O} \biggl(
\biggl(\log \frac{1}{x} \biggr)^{-1} \biggr) \biggr)
\]
as $x\rightarrow0$. Iterating this procedure $n-\bar{n}$ times using
the induction
argument, we finally get that
\[
\mathbb{P}\bigl[\mathrm{e}^{Y_1} + \cdots+ \mathrm{e}^{Y_{n}} \leq z\bigr] = \mathbb{P}\bigl[\mathrm{e}^{Y_1} + \cdots+ \mathrm{e}^{Y_{\bar{n}}} \leq z\bigr] \biggl(1+\mathrm{O} \biggl(
\biggl(\log \frac{1}{x} \biggr)^{-1} \biggr) \biggr)
\]
as $x\rightarrow0$, which completes the proof of the theorem, since
the asymptotics for
\[
\mathbb{P}\bigl[\mathrm{e}^{Y_1} + \cdots+ \mathrm{e}^{Y_{\bar{n}}} \leq z\bigr]
\]
can be computed using Lemma~\ref{easy.prop}.
\end{pf*}
\begin{pf*}{Proof of Corollary~\ref{laplace.cor}}
For any $u\in\mathbb{R}^n$, we have
%
\begin{equation}\label{probchange}
\mathbb E\bigl[\mathrm{e}^{\sum_{i=1}^n u_i Y_i}| X\leq x\bigr] = \frac{\mathbb
E[\mathrm{e}^{\sum_{i=1}^n u_i Y_i}\mathbf{1}_{X\leq x}]}{\mathbb{P}[X\leq x]} = \mathbb E
\bigl[\mathrm{e}^{\sum_{i=1}^n u_i Y_i}\bigr] \frac{\tilde{\mathbb{P}}[X\leq
x]}{\mathbb{P}[X\leq x]},
\end{equation}
where the symbol $\tilde{\mathbb{P}}$ stands for a new probability
determined from
\[
\frac{\mathrm{d}\tilde{\mathbb{P}}}{\mathrm{d}{\mathbb{P}}} = \frac{\mathrm{e}^{\sum_{i=1}^n u_i
Y_i}}{\mathbb E[\mathrm{e}^{\sum_{i=1}^n u_i Y_i}]}.
\]
Under this probability, $\ve Y \sim N(\veee{\mu}+ \cm\ve u, \cm)$. Applying
Theorem~\ref{gen.thm} to the numerator and the denominator of the
fraction in
\eqref{probchange}, and making cancellations, we get the following:
\begin{eqnarray*}
&& \mathbb E\bigl[\mathrm{e}^{\sum_{i=1}^n u_i Y_i}| X\leq x\bigr]
\\
&& \quad = x^{\sum_{i=1}^n u_i \sum_{j=1}^{\bar{n}}
\bar A_j b_{ij}} \mathrm{e}^{\veee{\mu}^\perp\ve u +
({1}/{2}) \ve u^\perp\cm\ve u} \frac{C_{\veee{\mu}+\cm\ve u,
\cm}}{C_{\veee{\mu},\cm}} \biggl(1+\mathrm{O} \biggl(
\biggl(\log \frac{1}{x} \biggr)^{-1} \biggr) \biggr)
\\
&&\quad = x^{\sum_{i=1}^n u_i \sum_{j=1}^{\bar{n}}
\bar A_j b_{ij}} \exp \Biggl\{\frac{1}{2}\sum
_{i,j=1}^n u_i u_j
\Biggl(b_{ij}-\sum_{
p,q=1}^{\bar{n}}
\bar a_{pq} b_{pi}b_{qj} \Biggr) \Biggr\}
\\
&& \qquad{}\times\exp \Biggl\{\sum_{i=1}^n
u_i \Biggl(\mu_i - \sum_{ p,q=1}^{\bar{n}}
\bar a_{pq}b_{pi} \biggl(\log\frac{\bar A_1+\cdots+\bar A_{\bar{n}}}{\bar A_q}+\bar
\mu_q \biggr) \Biggr) \Biggr\}
\\
&& \qquad{}\times\biggl(1+\mathrm{O} \biggl( \biggl(\log \frac{1}{x} \biggr)^{-1}
\biggr) \biggr)
\end{eqnarray*}
as $x\rightarrow0$. The symbols $C_{\veee{\mu}+\cm\ve u,\cm}$ and $C_{\veee{\mu},\cm}$, appearing in the previous estimates, stand for the constant
$C$ in Theorem~\ref{gen.thm},
evaluated for the log-normal variables, associated with the pairs
$(\veee{\mu}+\cm\ve u,\cm)$ and $(\veee{\mu},\cm)$, respectively.
It is easy to check that when $i\in\bar I$ or $j\in\bar I$,
necessarily
\[
b_{ij}=\sum_{ p,q=1}^{\bar{n}} \bar
a_{pq}b_{pi}b_{qj}.
\]
\upqed\end{pf*}
\begin{pf*}{Proof of Corollary~\ref{Ccorr}}
Recall the formula for the distribution function of $X$:
\begin{eqnarray*}
\mathbb{P}[X\leq x]  &=&  \int_{\mathrm{e}^{y_1}+\cdots+\mathrm{e}^{y_n}\leq x}\frac{1}{(2\uppi)^{{n}/{2}}\sqrt
{|\cm|}} \\
&&\hspace*{55pt}{}\times\exp
\biggl\{-\frac{1}{2} (\ve y-\veee{\mu})^\perp\cm^{-1} (
\ve y-\veee{\mu}) \biggr\} \,\mathrm{d}y_1\cdots \,\mathrm{d}y_n
\\
&=&\int_{\mathrm{e}^{z_1}+\cdots+\mathrm{e}^{z_n}\leq1} \frac{1}{(2\uppi)^{{n}/{2}}\sqrt{|\cm|}} \\
&&\hspace*{55pt}{}\times\exp \biggl\{-
\frac{1}{2} (\ve z + \mathbf{1} \log x-\veee{\mu})^\perp
\cm^{-1} (\ve z + \mathbf{1} \log x-\veee{\mu}) \biggr\}\, \mathrm{d}z_1
\cdots \,\mathrm{d}z_n.
\end{eqnarray*}
Differentiating with respect to $x$, we obtain an alternative
representation for the density:
\begin{eqnarray*}
p(x) &=&  \int_{\mathrm{e}^{z_1}+\cdots+\mathrm{e}^{z_n}\leq1} \frac{-\mathbf{1}^{\perp}\cm^{-1}(\ve z + \mathbf{1} \log x -
\veee{\mu})}{x(2\uppi)^{{n}/{2}}\sqrt{|\cm|}}
\\
&&\hspace*{35pt}\qquad {}\times \exp \biggl\{-\frac{1}{2} (\ve z + \mathbf{1} \log x-\veee{
\mu})^\perp\cm^{-1} (\ve z + \mathbf{1} \log x-\veee{\mu})
\biggr\}\, \mathrm{d}z_1\cdots \,\mathrm{d}z_n
\\
& =& -\frac{1}{x}\mathbb E\bigl[\mathbf{1}^\perp\cm^{-1}
(\ve Y-\veee{\mu}) \mathbf{1}_{X\leq
x}\bigr] = -\frac{1}{x}\mathbb E
\Biggl[\sum_{i=1}^n A_i
(Y_i - \mu_i)\mathbf 1_{X\leq
x} \Biggr].
\end{eqnarray*}
Next, we make a transformation inspired by Corollary~\ref{laplace.cor}.
Taking into account Remark~\ref{remgt1}, we see that
\begin{eqnarray*}
p(x) &=&  -\frac{1}{x}\mathbb E \Biggl[\sum_{i=1}^n
A_i \Biggl(Y_i - \log x \sum
_{j=1}^{\bar{n}} \bar A_j b_{ij} -
\mu_i \Biggr)\mathbf{1}_{X\leq
x} \Biggr]
\\
&&{}-\frac{\log x}{x}\mathbb E \Biggl[\sum_{i=1}^n
A_i \sum_{j=1}^{\bar{n}} \bar
A_j b_{ij} \mathbf{1}_{X\leq
x} \Biggr]
\\
& =& \frac{\mathbb{P}[X\leq x]}{x} \sum_{i=1}^n
A_i \mu_i - \frac{1}{x}\sum
_{i=1}^n A_i\mathbb E\bigl[
\bigl(Y_i - (1+\bar\lambda_i) \log x\bigr)
\mathbf{1}_{X\leq x}\bigr]
\\
&&{}-\frac{\log x}{x} \sum_{j=1}^{\bar{n}} \bar
A_j \mathbb{P}[X\leq x]
\\
& = &-\frac{\log x}{x} \sum_{j=1}^{\bar{n}}
\bar A_j \mathbb{P}[X\leq x] + \mathrm{O} \biggl(\frac{\mathbb{P}[X\leq x]}{x} \biggr)
\end{eqnarray*}
as $x\rightarrow0$. Here the constant $\bar\lambda_i$ is defined by
(\ref{Ela}).
In the reasoning above, we used the following estimate, which can be
derived from Corollary~\ref{laplace.cor}.
For every $i$, as $x\rightarrow0$,
\begin{eqnarray*}
&& \bigl\llvert \mathbb E\bigl[\bigl(Y_i - (1+\bar
\lambda_i) \log x\bigr) \mathbf{1}_{X\leq x}\bigr] \bigr\rrvert
\\
&& \quad \leq\mathbb{P}[X\leq x] \bigl(\mathbb E\bigl[\mathrm{e}^{Y_i - (1+\bar\lambda_i) \log x} | X\leq x\bigr] +
\mathbb E\bigl[\mathrm{e}^{-(Y_i - (1+\bar\lambda_i) \log x)} | X\leq x\bigr] \bigr)
\\
&& \quad  = C\mathbb{P}[X\leq x] \biggl(1+\mathrm{O} \biggl( \biggl(\log \frac{1}{x}
\biggr)^{-1} \biggr) \biggr)
\end{eqnarray*}
for some constant $C$.
\end{pf*}

\section{Asymptotic behavior of the right tail of a log-normal
difference}\label{Stail1}
In this section,\vspace*{1pt} we analyze the asymptotic behavior of the distribution
function and the density of the random variable $X^{(m)}$ as $x\to
+\infty$, assuming with no loss of
generality that $m\geq1$. If $m=0$, then the support of the distribution
of $X^{(m)}$ is $(-\infty,0)$, and the tail behavior at $0$ follows
from the results obtained in Section~\ref{Snachalo}.
\subsection{Right tail of a log-normal difference. Results and
discussions}\label{Srta}

Let us first consider, for every $p$ with $1\le p\le m$, the random variable
%
\begin{equation}\label{Esplit1}
X^{(1)}_p=\mathrm{e}^{Y_p}-\sum
_{k=m+1}^n\mathrm{e}^{Y_k}.
\end{equation}
Let $\Delta^p_{m,n}$ be the set of
weights ${\ve{w}}\in\mathbb{R}^n$ with $w_i=0$ for $ i=1,\ldots,m$,
$i\neq p$;
$w_p\geq0$; $w_i\leq0$ for $i=m+1,\ldots,n$; and $\sum w_i=1$. Let
$\bar{\ve{w}}_p\in\Delta^p_{m,n}$ be the unique point such that
%
\begin{equation}\label{minprobp}
\bar{\ve{w}}_p \cm\bar{\ve{w}}_p = \min
_{{\ve{w}}\in\Delta^p_{m,n}} {\ve{w}}^\perp\cm{\ve{w}},
\end{equation}
and define $\bar{n}^{(p)}$,
$\bar{I}^{(p)}$, $\bar k^{(p)}$, $\bar{\veee{\mu}}^{(p)}$, $\bar{\cm}^{(p)}$,
$\bar{a}_{ij}^{(p)}$, $\bar{A}^{(p)}$ as in equation \eqref{nbar.eq}
and below. We will say that assumption $(\mathcal A^p_1)$ holds if for
every $i\in\{m+1,\ldots,n\}\setminus\bar I^{(p)}$,
\[
\bigl(\ve{e}^i - \bar{\ve{w}}_p\bigr)^\perp\cm
\bar{\ve{w}}_p \neq0.
\]

It follows from Theorem~\ref{Tcharacter} (see Section~\ref{SStrap2}), that if assumption $(\mathcal A^p_1)$ is
satisfied, then
%
\begin{equation}\label{Esplitasymp}
\mathbb{P}\bigl[X^{(1)}_p\ge x\bigr]=\delta_{1,p}(
\log x)^{\delta_{2,p}}x^{\delta_{3,p}} \exp \biggl\{-\frac{\log^2x}{2 \delta_{4,p}} \biggr
\}\bigl(1+\mathrm{O} \bigl((\log x)^{-1} \bigr)\bigr)
\end{equation}
as $x\rightarrow\infty$, where
\begin{eqnarray*}
\delta_{1,p} &=& \frac{C^{(p)}}{\sum_{j=1}^{\bar{n}^{(p)}}\bar
{A}_j^{(p)}}, \qquad \delta_{2,p}=-
\frac{1+\bar{n}^{(p)}}{2},
\\
\delta_{3,p} &=& \sum_{i=1}^{\bar{n}_p}
\bar{A}_i^{(p)} \biggl(\log\frac{\sum_{j=1}^{\bar{n}^{(p)}}\bar{A}_j^{(p)}}{
|\bar{A}_i^{(p)}|}+\bar{
\mu}_i^{(p)} \biggr),\quad  \mbox{and}\quad \delta_{4,p}=
\Biggl(\sum_{j=1}^{\bar{n}^{(p)}}\bar{A}_j^{(p)}
\Biggr)^{-1}.
\end{eqnarray*}
%

In other words, the exponential rate of decay in the
leading term of the asymptotics of $\mathbb{P}[X_p^{(1)}\geq x]$ is
determined by the quantity
\[
\delta_{4,p} = \min_{\ve{w}\in\Delta^p_{m,n}}\ve{w}^\perp\cm
\ve{w}.
\]
Depending on the covariance matrix $\cm$, this rate may either be
equal to $b_{pp}$, the inverse of the variance of $Y_p$, in which case the
asymptotic behavior of $X^{(1)}_p$ is determined by $Y_p$ only, or be
greater than $b_{pp}$, in which case the asymptotic behavior of
$X^{(1)}_p$ is determined by more than one component of the vector
$(Y_1,\ldots,Y_n)$. One may therefore call the number $\delta_{4,p}$
the relative asymptotic variance of $Y_p$ with respect to
$Y_{m+1},\ldots,Y_n$.

The next assertion, which is one of the main results of the present
paper, shows that the asymptotic behavior of the distribution function
of the random variable $X^{(m)}$ is dominated by one (or several
similar) of the random variables $X^{(1)}_p$.
We will need the following parameters to describe the above-mentioned
domination:
%
\begin{eqnarray}
\label{P4}
\delta_4 &=& \max_{1\le p\le m}\delta_{4,p},
\qquad \mathcal{P}_4=\{p\dvt 1\le p\le m,\delta_{4,p}=
\delta_4\},
\\
\label{P3}
\delta_3 &=& \max_{p\in\mathcal{P}_4}\delta_{3,p},\qquad
\mathcal{P}_3=\{p\in\mathcal{P}_4\dvt \delta_{3,p}=
\delta_3\},
\\
\label{P2}
\delta_2 &=& \max_{p\in\mathcal{P}_3}\delta_{2,p},\qquad
\mathcal{P}_2=\{p\in\mathcal{P}_3\dvt \delta_{2,p}=
\delta_2\},
\end{eqnarray}
and finally
\[
\delta_1=\sum_{p\in\mathcal{P}_2}
\delta_{1,p}.
\]

\begin{theorem}\label{Tpresent}
Let assumption $(\mathcal A^p_1)$ hold for every
$p=1,\ldots,m$. Then the distribution function of the random variable
$X^{(m)}$ defined by (\ref{Ebass}) satisfies
%
\begin{equation}\label{dfdiff.eq}
\mathbb{P}\bigl[X^{(m)}\ge x\bigr]=\delta_1(\log
x)^{\delta_2}x^{\delta_3} \exp \biggl\{-\frac{\log^2x}{2 \delta_4} \biggr\}\bigl(1+\mathrm{O}
\bigl((\log x)^{-{1}/{2}} \bigr)\bigr)
\end{equation}
as $x\rightarrow\infty$.
\end{theorem}

\begin{remark}
When $m=n$, the variables $X^{(1)}_p$ are one-dimensional and log-normal,
and the result in Theorem~\ref{Tpresent} reduces to Theorem~1 of
\cite{AR-N} which shows that the asymptotic
behavior of the right tail of $\mathrm{e}^{Y_1}+\cdots+\mathrm{e}^{Y_n}$ is determined by
the components of $(Y_1,\ldots,Y_n)$ which have the largest variance.
For other values of $m$, Theorem~\ref{Tpresent} extends Theorem~1 of
\cite{AR-N} by showing that the asymptotic behavior of the right tail
of $X^{(m)}$
is determined by the components of $(Y_1,\ldots,Y_m)$, which have the
largest relative asymptotic variance with respect to
$(Y_{m+1},\ldots,Y_n)$.
\end{remark}

As in the case of Theorem~\ref{gen.thm}, several useful corollaries can
be derived from Theorem~\ref{Tpresent}. We omit the proofs of those
corollaries, since they are very
similar to those given in Section~\ref{SSfmr1}.

\begin{corollary}\label{laplace2.cor}
Suppose that assumption\ $(\mathcal A^p_1)$ holds for every $p=1,\ldots,m$, and that
the set $\mathcal P_4$ is a singleton, $\mathcal P_4 = \{p\}$. Then, as
$x\to\infty$, for any $\ve u\in\mathbb{R}^n$,
%
\begin{eqnarray}
&& \mathbb E \bigl[\mathrm{e}^{\sum_{i=1}^n u_i Y_i}| X\geq x \bigr]\nonumber\\
&&\quad =  \mathbb E
\bigl[\mathrm{e}^{\sum_{i=1}^n u_i Y_i}| X= x \bigr]\nonumber\\
&&\quad  = x^{\sum_{j=1}^n u_j \sum_{i=1}^{\bar{n}^{(p)}}\bar
A_i^{(p)}b_{\bar k^{(p)}(i),j} }\nonumber
\\[-8pt]
\label{lt2.eq}\\[-8pt]
\nonumber
&&\qquad {}\times\exp \Biggl\{ \sum_{j=1}^n
u_j \Biggl(\mu_j-\sum_{i,k=1}^{\bar{n}^{(p)}}
\bar a^{(p)}_{ik} b_{\bar k^{(p)}(i),j} \biggl(\log
\frac{\sum_{j=1}^{\bar{n}^{(p)}} \bar A^{(p)}_j}{|\bar A^{(p)}_k|}+ \mu_{\bar k^{(p)}(k)} \biggr) \Biggr) \Biggr\}
\\
&&\qquad {}\times\exp \Biggl\{+\frac{1}{2} \sum_{j,l\notin\bar I^{(p)}}^n
u_j u_l \Biggl(b_{jl}-\sum
_{i,k=1}^{\bar{n}^{(p)}} \bar a^{(p)}_{ik}
b_{\bar k^{(p)}(i),j} b_{\bar
k^{(p)}(k),l} \Biggr) \Biggr\} \nonumber
\\
&&\qquad {}\times \biggl(1+\mathrm{O} \biggl( \biggl(\log \frac{1}{x} \biggr)^{-1}
\biggr) \biggr).\nonumber
\end{eqnarray}
\end{corollary}

In the next statement, we use the notation introduced before the
formulation of Theorem~\ref{Tpresent}.

\begin{corollary}\label{Ccorrr}
Suppose assumption $(\mathcal A^p_1)$ holds for every $p=1,\ldots,m$.
Then, as $x\to\infty$, the density $p^{(m)}$ of the random variable
$X^{(m)}$ satisfies
%
\begin{equation}\label{densdiff.eq}
p^{(m)}(x) = \frac{\delta_1}{\delta_4}(\log x)^{\delta_2+1}x^{\delta_3-1}
\exp \biggl\{-\frac{\log^2 x}{2\delta_4} \biggr\}\bigl(1+\mathrm{O} \bigl((\log x)^{-{1}/{2}}
\bigr)\bigr).
\end{equation}
\end{corollary}

\subsection{Right tail of a log-normal difference. Proofs}\label{SStrap2}
Let us first consider the case when $m=1$. Define
%
\begin{equation}\label{newdeltan}
\Delta_{1,n}:= \Biggl\{{\ve{w}}\in\mathbb{R}^n\dvt
w_1\geq0, w_i \leq0, i=2,\ldots,n, \mbox{and } \sum
_{i=1}^n w_i = 1\Biggr\},
\end{equation}
and introduce $\bar{\ve{w}}\in\Delta_{1,n}$ as the unique
vector such that
\[
\bar{\ve{w}}^\perp\cm\bar{\ve{w}} = \min_{{\ve{w}}\in\Delta_{1,n}} {
\ve{w}}^\perp\cm {\ve{w}}.
\]
The existence and uniqueness
of $\bar{\ve{w}}$ follows from the non-degeneracy of $\cm$.
When $A_1>0$
and $A_k<0$ for $k=2,\ldots,n$, $\bar{\ve{w}}$ is given by
\eqref{easycase.eq}.
In the general case, we define $\bar{n}$, $\bar I$, $\bar
{\veee{\mu}}$, $\bar\cm$, $\bar a_{ij}$ and $\bar A$ as in equation \eqref
{nbar.eq} and below.

Since by the definition of $\Delta_{1,n}$, $\bar{w}_1>0$, and moreover the
variables $Y_2,\ldots,Y_n$ are exchangeable in the definition of
$X^{(1)}$, we shall assume with no loss of generality that $\bar I = \{
1,\ldots,\bar{n}\}$.
This has already been done in Section~\ref{SSfmr1}.

Observe that the minimum in
$ {\min_{{\ve{w}}\in\Delta_{1,\bar{n}}} {\ve{w}}^\perp\bar
\cm{\ve{w}}}$
is attained in the
interior of $\mathbb{R}_+^{\bar{n}}$. This implies that $\bar A_1>0$ and
$\bar A_k<0$ for $k=2,\ldots,\bar{n}$ (see a similar reasoning in
Section~\ref{SSfmr1}).

The following preliminary lemma concerns the case where $\bar{n} = n$.

\begin{lemma}\label{Tstandard}
Let $A_1> 0$, $A_2< 0,\ldots, A_n< 0$. Then the following formulas hold:
%
\begin{eqnarray}
p^{(1)}(x)&= & C(\log x)^{({1-n})/{2}} x^{-1+\sum_{i=1}^n A_i ( \log{\sum_{j=1}^n A_j}/{|A_i|} +
\mu_i )}
\nonumber
\\[-8pt]
\label{Efinf}\\[-8pt]
\nonumber
&& {}\times\exp \Biggl\{-\frac{1}{2}\log^2x\sum
_{j=1}^nA_j \Biggr\} \bigl(1+\mathrm{O} \bigl(
(\log x )^{-1} \bigr) \bigr),
\\
\mathbb{P}\bigl[X^{(1)}\geq x\bigr] &=&  \frac{C}{A_1+\cdots+A_n}(\log
x)^{-({1+n})/{2}} x^{\sum_{i=1}^n A_i ( \log{\sum_{j=1}^n A_j}/{|A_i|} +
\mu_i )}
\nonumber
\\[-8pt]
\\[-8pt]
\nonumber
&&{}\times\exp\Biggl\{-\frac{1}{2}\bigl(\log^2x\bigr)\sum
_{j=1}^nA_j\Biggr\} \bigl(1+\mathrm{O}
\bigl( (\log x )^{-1} \bigr) \bigr)
\end{eqnarray}
as $x\rightarrow\infty$. The constant $C$ in (\ref{Efinf}) is given by
\begin{eqnarray*}
C &=& \exp \Biggl\{-\frac{1}{2}\sum_{i,m=1}^na_{im}
\biggl( \log\frac{\sum_{j=1}^n A_j}{|A_i|} + \mu_i \biggr) \biggl( \log
\frac{\sum_{j=1}^n A_j}{|A_m|} + \mu_m \biggr) \Biggr\}
\\
&&{}\times\frac{1}{\sqrt{2\uppi|\cm|}} \sqrt{\frac{\sum_{j=1}^n A_j}{\prod_{i=1}^n |A_i|}}.
\end{eqnarray*}
\end{lemma}

\begin{pf}
Differentiating the distribution function, we obtain the following
representation of the
density $p^{(1)}$ of $X^{(1)}$:
%
\begin{eqnarray}
p^{(1)}(x) &=& \frac{1}{(2\uppi)^{{n}/{2}}\sqrt{|\cm|}} \exp \Biggl\{-\frac{1}{2}
\log^2x\sum_{j=1}^nA_j
\Biggr\}
\nonumber
\\[-8pt]
\label{Ef4}
\\[-8pt]
\nonumber
&& {}\times\int_{D^{1,n-1}_1}\widetilde{\Phi}(x_1,
\ldots,x_{n-1}) \exp \bigl\{-\log x \widetilde{\Psi}(x_1,
\ldots,x_{n-1}) \bigr\} \,\mathrm{d}x_1\cdots \,\mathrm{d}x_{n-1}
\end{eqnarray}
with
\begin{eqnarray*}
\widetilde{\Phi}(x_1,\ldots,x_{n-1})&=&\Phi(x_1,
\ldots ,x_{n-1},x_1-x_2-\cdots-x_{n-1}-1),
\\
\widetilde{\Psi}(x_1,\ldots,x_{n-1}) &=& \Psi(x_1,
\ldots ,x_{n-1},x_1-x_2-\cdots-x_{n-1}-1),
\\
\Phi(x_1,\ldots,x_n) &=& \frac{1}{x_1\cdots x_n}\exp \Biggl\{-
\frac{1}{2}\sum_{i,j=1}^na_{ij}(
\log x_i-\mu_i) (\log x_j-\mu_j)
\Biggr\},
\\
\Psi(x_1,\ldots,x_n) &=& \sum_{i=1}^nA_i(
\log x_i-\mu_i)
\end{eqnarray*}
and
\begin{eqnarray*}
D^{1,n-1}_1=\bigl\{(x_1,\ldots,x_{n-1})
\in\mathbb{R}^{n-1} \dvt x_i\ge0, 1\le i\le n-1;
x_1-x_{2}-\cdots-x_{n-1}> 1\bigr\}.
\end{eqnarray*}
We have for all $1\le i\le n-1$,
%
\begin{equation}\label{Eder}
\frac{\partial\widetilde\Psi}{\partial x_i}=\frac{A_i}{x_i}+\frac
{A_n}{x_1-x_2-\cdots- x_{n-1}-1}s_i,
\end{equation}
where $s_1=1$ and $s_i=-1$ for $1< i\le n-1$. Now, it is easy to see
that the solution $\ve x^{*}$ to
the system of equations $\frac{\partial\widetilde\Psi}{\partial
x_i}=0$, $1\le i\le n-1$, is given by
%
\begin{equation}\label{Ef7}
x^{*}_i=\frac{A_i}{\sum_{j=1}^nA_j}s_i, \qquad 1\le i\le
n-1.
\end{equation}

Under\vspace*{1pt} the assumptions in the formulation of the lemma, it is clear that
$\ve x^{*}$ belongs to
the interior of the set $D^{1,n-1}_1$. We will next apply
Laplace's method to the integral in (\ref{Ef4}). However, first we
need to check that the Hessian matrix
of the function $\widetilde{\Psi}$ at the point $\ve x^*$, that is, the matrix
$H(\ve x^{*}):=[h_{im}]_{i,m=1,\ldots, n-1}$, is positive-definite.
It is not hard to see that
\[
h_{im}=-\frac{ (\sum_{j=1}^nA_j )^2}{A_n}\qquad  \mbox{if } i\neq m,
\]
and
\[
h_{ii}=- \Biggl(\sum_{j=1}^nA_j
\Biggr)^2 \biggl(\frac{1}{A_i}+\frac
{1}{A_n} \biggr).
\]
Therefore,
%
\begin{equation}\label{Eeco}
H\bigl(\ve x^{*}\bigr)=-\frac{ (\sum_{j=1}^nA_j )^2}{A_n}J,
\end{equation}
where $J$ is the $(n-1)\times(n-1)$-matrix with the entries $1+\frac
{A_n}{A_1},\ldots,1+\frac{A_n}{A_{n-1}}$
along the main diagonal, and all the entries outside the main diagonal
equal to 1. It is an easy exercise in linear algebra to show that the
leading principal minor of order $p$ of the matrix $J$
is equal to
%
\begin{equation}\label{Eg}
A_n^{p-1} \Biggl(\prod_{i=1}^p
\frac{1}{A_i} \Biggr) \Biggl(\sum_{m=1}^pA_m+A_n
\Biggr).
\end{equation}
Recall that it is assumed in Lemma~\ref{Tstandard} that
%
\begin{equation}\label{Eco}
A_1> 0,\qquad  A_2< 0,\qquad  \ldots,\qquad A_n< 0.
\end{equation}
Under this assumption, the numbers in (\ref{Eg}) are positive for
$p=1,\ldots,n-1$. For instance, if $p=1$, then we have
\[
\frac{1}{A_1}(A_1+A_n)=\frac{1}{A_1} \Biggl[
\sum_{m=1}^nA_m-\sum
_{m=2}^{n-1}A_m \Biggr]> 0.
\]
The previous inequality follows from the positive-definiteness of the
matrix $\cm^{-1}$ and condition (\ref{Eco}).
It follows that the matrix $J$ is positive definite, and hence the
matrix $H(\ve x^*)$ is also positive
definite. Moreover, the determinant of $H(\ve x^*)$ is given by
\[
\frac{ (\sum_{j=1}^n A_j )^{2n-1}}{\prod_{i=1}^n |A_i|}.
\]
Next, taking in the account what was said above, we see that Laplace's
method can be applied to the integral in
(\ref{Ef4}). Similarly to \eqref{Evis6}, we get the following formula:
\begin{eqnarray*}
p^{(1)}(x)&=&\frac{1}{(2\uppi)^{{n}/{2}}x\sqrt{|\cm|}} \exp \Biggl\{-\frac{1}{2}
\log^2x\sum_{j=1}^nA_j
\Biggr\} \frac{1}{\sqrt{\det(H({\ve x}^{*}))}} \biggl(\frac{2\uppi}{\log x} \biggr)^{({n-1})/{2}}
\\
&&{}\times\widetilde\Phi\bigl(x_1^{*},\ldots,x_{n-1}^{*}
\bigr) \exp \bigl\{-\log x \widetilde\Psi\bigl(x_1^{*},
\ldots,x_{n-1}^{*}\bigr) \bigr\} \bigl(1+\mathrm{O} \bigl( (\log x
)^{-1} \bigr) \bigr)
\\
& = & (\log x)^{({1-n})/{2}} x^{-1+\sum_{i=1}^n A_i ( \log{\sum_{j=1}^n A_j}/{|A_i|} + \mu
_i )}\exp \Biggl\{-\frac{1}{2}
\log^2x\sum_{j=1}^nA_j
\Biggr\}
\\
&&{}\times\exp \Biggl\{-\frac{1}{2}\sum_{i,m=1}^na_{im}
\biggl( \log\frac
{\sum_{j=1}^n A_j}{|A_i|} + \mu_i \biggr) \biggl( \log
\frac{\sum_{j=1}^n A_j}{|A_m|} + \mu_m \biggr) \Biggr\}
\\
&&{}\times\frac{1}{\sqrt{2\uppi|\cm|}} \sqrt{\frac{\sum_{j=1}^n A_j}{\prod_{i=1}^n |A_i|}} \bigl(1+\mathrm{O} \bigl(
(\log x )^{-1} \bigr) \bigr)
\end{eqnarray*}
as $x\rightarrow\infty$. The asymptotic behavior of the distribution
function can be characterized by
integrating the asymptotic formula for the density.
\end{pf}

We will next focus on the case where $m$ is still equal to one, but the equality
$\bar{n} = n$ may not hold. Our next result requires the following assumption:
\begin{enumerate}[$(\mathcal A_1)$]
\item[$(\mathcal A_1)$] For every $i\in\{1,\ldots,n\}\setminus\bar I$,
\[
\bigl(\ve{e}^i - \bar{\ve{w}}\bigr)^\perp\cm\bar{\ve{w}}
\neq0.
\]
\end{enumerate}

\begin{remark}\label{rema1}
Assumption $(\mathcal A_1)$, although it has the same form as
assumption $(\mathcal
A)$ above, is a different assumption on the covariance matrix $\cm$,
because the weight vector $\bar{\ve{w}}$ is computed differently now
(with $\Delta_{1,n}$ instead of $\Delta$).
Assumption $(\mathcal A_1)$ is equivalent to the following:
For every $i\in\{1,\ldots,n\}\setminus\bar I$,
%
\begin{equation}\label{Eopp}
\bigl(\ve{e}^i - \bar{\ve{w}}\bigr)^\perp\cm\bar{\ve{w}} <
0.
\end{equation}
Indeed, since $1\in\bar{I}$, for every $i\in\{1,\ldots,n\}\setminus
\bar I$ and for $\varepsilon>0$ small enough, $\bar{\ve{w}} -\varepsilon
(\ve{e}^i -
\bar{\ve{w}})$ belongs to $\Delta_{1,n}$. Therefore, the inequality
opposite to that in (\ref{Eopp}) would
lead to a contradiction to the fact that $\bar{\ve{w}}$ is a minimizer.
\end{remark}

\begin{theorem}\label{Tcharacter}
Let assumption $(\mathcal A_1)$ hold true. Then,
as $x\to+\infty$,
%
\begin{eqnarray}
\mathbb{P}\bigl[X^{(1)}\geq x\bigr] &=& \frac{C}{\bar A_1+\cdots+\bar A_{\bar{n}}}(\log
x)^{-({1+{\bar{n}}})/{2}} x^{\sum_{i=1}^{\bar{n}} \bar A_i ( \log{\sum_{j=1}^{\bar{n}}
\bar A_j}/{|\bar A_i|} +
\bar\mu_i )}
\nonumber
\\[-8pt]
\\[-8pt]
\nonumber
&&{}\times\exp\Biggl\{-\frac{1}{2}\bigl(\log^2x\bigr)\sum
_{j=1}^{\bar{n}}\bar A_j\Biggr\}
\bigl(1+\mathrm{O} \bigl( (\log x )^{-1} \bigr) \bigr),
\end{eqnarray}
where
\[
C=\exp \Biggl\{-\frac{1}{2}\sum_{i,m=1}^{\bar{n}}
\bar a_{im} \biggl( \log \frac{\sum_{j=1}^{\bar{n}} \bar A_j}{|\bar A_i|} + \bar\mu_i
\biggr) \biggl( \log\frac{\sum_{j=1}^{\bar{n}} \bar A_j}{|\bar A_m|} + \bar\mu _m \biggr) \Biggr\}
\times\frac{1}{\sqrt{2\uppi|\bar\cm|}} \sqrt{\frac{\sum_{j=1}^{\bar{n}}
\bar A_j}{\prod_{i=1}^{\bar{n}} |\bar A_i|}}.
\]
\end{theorem}

\begin{pf}
We will only sketch the proof, which is very similar to that of Theorem~\ref{gen.thm}.
If $\bar{n}=n$, the result follows from Lemma~\ref{Tstandard}. Next,
assume that
$k \in\{\bar{n}, \ldots, n-1\}$, $x>1$, and
let $a,b$ be such that $x = \mathrm{e}^a - \mathrm{e}^b$. On the one hand, clearly,
\[
\mathbb{P}\bigl[\mathrm{e}^{Y_1}\geq \mathrm{e}^{Y_2} +\cdots+ \mathrm{e}^{Y_n}+x
\bigr] \leq\mathbb{P}\bigl[\mathrm{e}^{Y_1}\geq \mathrm{e}^{Y_2} +\cdots+
\mathrm{e}^{Y_{n-1}}+x\bigr].
\]
On the other hand,
\begin{eqnarray*}
&& \mathbb{P}\bigl[\mathrm{e}^{Y_1}\geq \mathrm{e}^{Y_2} +\cdots+ \mathrm{e}^{Y_{k+1}}+x
\bigr] \\
&& \qquad \geq  \mathbb{P}\bigl[\mathrm{e}^{Y_2} + \cdots+ \mathrm{e}^{Y_{k}} \leq
\mathrm{e}^{Y_1} - \mathrm{e}^a, Y_{k+1} \leq b\bigr]
\\
&& \qquad =  \mathbb{P}\bigl[\mathrm{e}^{Y_2} + \cdots+ \mathrm{e}^{Y_{k}}\leq
\mathrm{e}^{Y_1} - \mathrm{e}^a\bigr] - \mathbb{P}\bigl[\mathrm{e}^{Y_2} +
\cdots+ \mathrm{e}^{Y_{k}}\leq \mathrm{e}^{Y_1} -\mathrm{e}^a, Y_{k+1} >
b\bigr].
\end{eqnarray*}
Since $k\geq\bar{n}$, $\sum_{i=1}^k \bar{w}_i = 1$. Moreover, since
$\bar{w}_1>0$, we may define $\tilde w_1 = \frac{1}{w_1}$ and $\tilde
w_i = -\frac{w_i}{w_1}$ for $i=1,\ldots,k$; these weights are positive
and satisfy $\sum_{i=1}^k \tilde w_i = 1$.

We have
\begin{eqnarray*}
&& \mathbb{P}\bigl[\mathrm{e}^{Y_2} + \cdots+ \mathrm{e}^{Y_{k}}\leq \mathrm{e}^{Y_1}
- \mathrm{e}^a, Y_{k+1} > b\bigr]
\\
&& \quad \leq\mathbb{P}\bigl[\tilde w_1 \mathrm{e}^a + \tilde
w_2 \mathrm{e}^{Y_2} + \cdots+ \tilde w_k
\mathrm{e}^{Y_{k}}\leq \mathrm{e}^{Y_1}, Y_{k+1} > b\bigr]
\\
&&\quad  \leq\mathbb{P}[\tilde w_1 a + \tilde w_2
{Y_2} + \cdots+ \tilde w_k {Y_{k}}
\leq{Y_1}, Y_{k+1} > b]
\\
&& \quad = \mathbb{P}[a \leq \bar{w}_1 Y_1 + \bar{w}_2
{Y_2} + \cdots+ \bar{w}_k {Y_{k}},
Y_{k+1} > b]
\\
&& \quad \leq\mathbb{P}\Biggl[a -\sum_{i=1}^k
\bar{w}_i Y_i \leq\alpha( Y_{k+1} - b)\Biggr]
\\
&& \quad  = \mathbb{P}\Biggl[\sum_{i=1}^k \bar{w}_i Y_i +\alpha Y_{k+1} \geq a +\alpha b \Biggr]
\\
&& \quad  = N \biggl(\frac{-a-\alpha b + \mathbb E[\bar{\ve{w}}^\perp\ve Y
+ \alpha Y_{k+1}]}{\sqrt{\operatorname{Var} [\bar{\ve{w}}^\perp\ve Y + \alpha
Y_{k+1}]}} \biggr)
\end{eqnarray*}
for all $\alpha>0$.

Next, reasoning as in the proof of Theorem~\ref{gen.thm}, let
\[
x_{k+1} = \frac{({\ve e}^{k+1})^\perp\cm\bar{\ve{w}}}{\bar{\ve{w}}^\perp\cm
\bar{\ve{w}}} < 1.
\]
Then,
\[
N \biggl(\frac{-a-\alpha b + \mathbb E[\bar{\ve{w}}^\perp\ve Y + \alpha
Y_{k+1}]}{\sqrt{\operatorname{Var} [\bar{\ve{w}}^\perp\ve Y + \alpha
Y_{k+1}]}} \biggr)= N \biggl(\frac{-a-\alpha b + \bar
{\ve{w}}^\perp\veee{\mu}+ \alpha\mu_{k+1}}{\sqrt{\bar{\ve{w}}^\perp\cm
\bar{\ve{w}}  (1+2
\alpha x_{k+1} +{\alpha^2 \cm_{k+1,k+1}}/({\bar{\ve{w}}^\perp\cm
\bar{\ve{w}}}))}} \biggr).
\]
Now, we take
\[
b =\frac{1}{2} (x_{k+1}+1) \log x \quad  \Rightarrow \quad  a = \log
\bigl(x+\mathrm{e}^b\bigr) = \log x + \log \bigl(1+ x^{-({1}/{2})(1-x_{k+1})} \bigr).
\]
Using these substitutions, we obtain
\begin{eqnarray*}
&& \frac{-a-\alpha b + \bar
{\ve{w}}^\perp\veee{\mu}+ \alpha\mu_{k+1}}{\sqrt{\bar{\ve{w}}^\perp\cm
\bar{\ve{w}}  (1+2
\alpha x_{k+1} +{\alpha^2 \cm_{k+1,k+1}}/({\bar{\ve{w}}^\perp
\cm\bar{\ve{w}}}) )}}
\\
&& \quad = \frac{-\log x(1+({\alpha}/{2})(1+x_{k+1}))- \log (1+ x^{-({1}/{2})(1-x_{k+1})} )+ \bar
{\ve{w}}^\perp\veee{\mu}+ \alpha\mu_{k+1}}{\sqrt{\bar{\ve{w}}^\perp\cm
\bar{\ve{w}}  (1+2
\alpha x_{k+1} +{\alpha^2 \cm_{k+1,k+1}}/({\bar{\ve{w}}^\perp
\cm\bar{\ve{w}}}) )}}
\\
&& \quad \leq-\frac{\log x}{\sqrt{\bar{\ve{w}}^\perp\cm\bar{\ve{w}}}} \frac
{1+({\alpha}/{2})(1+x_{k+1})}{\sqrt{1+2
\alpha x_{k+1} +{\alpha^2 \cm_{k+1,k+1}}/({\bar{\ve{w}}^\perp
\cm\bar{\ve{w}}})}} + C_{k+1},
\end{eqnarray*}
where $C_{k+1}$ is a constant independent of $x$. Next, reasoning as in the
proof of Theorem~\ref{gen.thm}, we see that there exist
$\alpha$ small enough and $\varepsilon_{k+1}>0$ such that for all $x>1$,
\[
\frac{-a-\alpha b + \mathbb E[\bar{\ve{w}}^\perp\ve Y + \alpha
Y_{k+1}]}{\sqrt{\operatorname{Var} [\bar{\ve{w}}^\perp\ve Y + \alpha
Y_{k+1}]}} \leq-\frac{\log x}{\sqrt{\bar{\ve{w}}^\perp\cm\bar{\ve{w}}}}(1+\varepsilon_{k+1})+C_{k+1}.
\]
The rest of the proof is similar to that of Theorem~\ref{gen.thm}
modulo some trivial changes.
\end{pf}

\begin{pf*}{Proof of Theorem~\ref{Tpresent}}
It is clear that the sets $\mathcal{P}_4$, $\mathcal{P}_3$, and $\mathcal{P}_2$
defined by (\ref{P4}),
(\ref{P3}), and (\ref{P2}), respectively, are not empty.

\textit{Upper estimate}.  Fix a positive function $\varphi$ such that
$\varphi(x)\to
0$ as $x\to\infty$. Then we have
%
\begin{eqnarray}
\mathbb{P}\bigl[X^{(m)}\geq x\bigr]  &\leq & \sum
_{1\leq p\leq m} \mathbb{P}\bigl[X_{p} \geq \bigl(1-
\varphi(x)\bigr) (X_{m+1}+\cdots+X_n + x)\bigr]\notag
\nonumber\\
\label{upperXm} &&{}+ \sum_{1\leq p,q \leq m, p\neq q} \mathbb{P} \biggl[X_p
\geq\frac
{\varphi(x)}{m-1}(X_{m+1}+\cdots+X_n+x),
\\
&&\qquad \hspace*{57pt} X_q \geq \frac{\varphi(x)}{m-1}(X_{m+1}+
\cdots+X_n+x) \biggr].\nonumber
\end{eqnarray}
Formula (\ref{upperXm}) can be established as follows. Let $E_1$,
$E_2$, and $F_i$ with $1\le i\le m$, be random variables. Then it is
not hard to prove that the following set theoretical inclusion holds:
%
\begin{eqnarray}
\bigl\{F_1+\cdots+F_m\ge E_1+(m-1)E_2
\bigr\} &\subset & \bigcup_{p=1}^m
\{F_p\ge E_1\}
\nonumber
\\[-8pt]
\label{Ekr1}
\\[-8pt]
\nonumber
&&{}\cup  \biggl[\bigcup_{1\leq p,q \leq m, p\neq q}
\{F_p\ge E_2,F_q\ge E_2\}
\biggr].
\end{eqnarray}
Next, using (\ref{Ekr1}) with
\[
F_p=X_p, \qquad  E_1=\bigl(1-\varphi(x)\bigr)
(X_{m+1}+\cdots+X_n + x),
\]
and
\[
E_2=\frac{\varphi(x)}{m-1}(X_{m+1}+\cdots+X_n+x),
\]
and taking into account the countable subadditivity of $\mathbb{P}$, we
obtain (\ref{upperXm}).

To estimate the terms in the first sum in formula (\ref{upperXm}), we
introduce the following probability measure:
%
\begin{equation}\label{Esu}
\frac{\mathrm{d}\widetilde{\mathbb{P}}}{\mathrm{d}\mathbb{P}} = \frac{\mathrm{e}^{\log(1-\varphi(x))
\ve{e}^p
\cm^{-1}\ve{Y}}}{\mathbb E[\mathrm{e}^{\log(1-\varphi(x)) \ve{e}^p
\cm^{-1}\ve{Y}}]},
\end{equation}
where ${\ve Y}=(Y_1,\ldots,Y_n)$, and $\ve{e}^p$ is the vector with
$p$th component equal to $1$ and all the
other components equal to zero.
Note that the measure $\widetilde{\mathbb{P}}$ depends on $p$. However,
we omit the parameter $p$ in the symbol
$\widetilde{\mathbb{P}}$ for the sake of simplicity. It is not hard to
see that under the probability~$\widetilde{\mathbb{P}}$, we have
$\ve{Y}\sim N(\veee{\mu}+\log(1-\varphi(x))\ve{e}^p,\cm)$. In other words,
the law of the random vector
\[
\bigl(Y_p -\log\bigl(1-\phi(x)\bigr), Y_{m+1},
\ldots,Y_n \bigr)
\]
under $\widetilde{\mathbb{P}}$ coincides with the law of the random
vector $(Y_p , Y_{m+1},\ldots,Y_n )$ under ${\mathbb{P}}$,
which means that
%
\begin{equation}\label{nophi.eq}
\widetilde{\mathbb{P}}\bigl[X_{p} \geq \bigl(1-\varphi(x)\bigr)
(X_{m+1}+\cdots+X_n + x)\bigr] = {\mathbb{P}}[X_{p} \geq X_{m+1}+\cdots+X_n +
x].
\end{equation}
It follows from (\ref{Esu}) that
%
\begin{eqnarray}
&& \mathbb{P}\bigl[X_{p} \geq \bigl(1-\varphi(x)\bigr)
(X_{m+1}+\cdots+X_n + x)\bigr]
\nonumber
\\
&&  \label{Ehr}\quad = \widetilde{\mathbb E} \biggl[\frac{\mathrm{d}{\mathbb{P}}}{\mathrm{d}\widetilde{\mathbb{P}}} \mathbf{1}_{\{X_{p}
\geq
(1-\varphi(x))(X_{m+1}+\cdots+X_n + x)\}}
\biggr]
\\
&&\quad  = \mathbb E\bigl[\mathrm{e}^{\log(1-\varphi(x)) \ve{e}^{p}
\cm^{-1}\ve{Y}}\bigr] \widetilde{\mathbb E} \bigl[
\mathrm{e}^{-\log(1-\varphi(x))\ve{e}^{p} \cm^{-1}\ve{Y}}\mathbf{1}_{\{
X_{p} \geq
(1-\varphi(x))(X_{m+1}+\cdots+X_n + x)\}} \bigr].\nonumber
\end{eqnarray}
Next, let $r$ and $q$ be positive numbers satisfying
$\frac{1}{r}+\frac{1}{q} = 1$. Then, using (\ref{Ehr}), \eqref{nophi.eq} and H\"{o}lder's inequality, we obtain
\begin{eqnarray*}
&& \mathbb{P}\bigl[X_{p} \geq \bigl(1-\varphi(x)\bigr) (X_{m+1}+
\cdots+X_n + x)\bigr]
\\
&& \quad \leq\mathbb E\bigl[\mathrm{e}^{\log(1-\varphi(x))\ve{e}^{p}
\cm^{-1}\ve{Y}}\bigr] \widetilde{\mathbb E} \bigl[
\mathrm{e}^{-r\log(1-\varphi(x)) \ve{e}^{p}
\cm^{-1}\ve{Y}} \bigr]^{1/r}
\\
&&\qquad  {}\times\widetilde{\mathbb{P}}\bigl[X_{p} \geq \bigl(1-\varphi(x)
\bigr) (X_{m+1}+\cdots+X_n + x)\bigr]^{1/q}
\\
&& \quad  = \mathbb E\bigl[\mathrm{e}^{\log(1-\varphi(x)) \ve{e}^{p}
\cm^{-1}\ve{Y}}\bigr]^{1-{1}/{r}} {\mathbb E} \bigl[
\mathrm{e}^{-(r-1)\log(1-\varphi(x)) \ve{e}^{p}
\cm^{-1}\ve{Y}} \bigr]^{1/r}
\\
&&\qquad  {}\times{\mathbb{P}}[X_{p} \geq X_{m+1}+
\cdots+X_n + x]^{1/q}
\\
&&\quad  = \mathrm{e}^{({r-1})/{2} \log^2(1-\varphi(x))a_{pp}}{\mathbb{P}}\bigl[X_{p} \geq (X_{m+1}+
\cdots+X_n + x)\bigr]^{1/q}.
\end{eqnarray*}
Next, set
%
\begin{equation}\label{Evar}
\varphi(x) = \frac{1}{\log^2 x}, \qquad  r = r(x) = \log^3 x,\quad  \mbox{and}\quad
\frac{1}{q(x)} =1-\frac{1}{\log^3 x}.
\end{equation}
Then we have
\begin{eqnarray*}
\exp \biggl\{\frac{r(x)-1}{2} \log^2\bigl(1-\varphi(x)
\bigr)a_{pp} \biggr\}& =& \exp \biggl\{\frac{\log^3(x)-1}{2}
\log^2\bigl(1-\log^{-2}(x)\bigr)a_{pp} \biggr\}
\\
& =& 1+\mathrm{O} \biggl(\frac{1}{\log x} \biggr),
\end{eqnarray*}
and hence, by Theorem~\ref{Tcharacter},
\begin{eqnarray*}
&& \mathbb{P} \bigl[X_{p} \geq\bigl(1-\varphi(x)\bigr)
(X_{m+1}+\cdots+X_n +x) \bigr]
\\
&& \quad \leq\delta _{1,p}(\log x)^{\delta_{2,p}}x^{\delta_{3,p}} \exp \biggl
\{-\frac{\log^2x}{2 \delta_{4,p}} \biggr\}\bigl(1+\mathrm{O} \bigl((\log x)^{-1}\bigr) \bigr)
\end{eqnarray*}
as $x\rightarrow\infty$.

It remains to estimate the terms in the second sum in
\eqref{upperXm}.
For any two integers $p$ and $q$ with $1\leq p,q\leq m$ and $p\neq q$,
let $\Delta^{p,q}_{m,n} $ be the set of weights $\ve{w}\in\mathbb{R}^n$
with $w_i = 0$ for $i=1,\ldots,m,
i\neq p$, $i\neq q$; $w_p \geq0$, $w_q\geq0$; $w_i \leq0$ for
$i=m+1,\ldots,n$; and $\sum
w_i =1$. Recall that
\[
\Delta^{p}_{m,n} = \bigl\{\ve{w}\in\Delta^{p,q}_{m,n}\dvt w_q = 0\bigr\}.
\]
By Jensen's inequality, for any ${\ve{w}}\in\Delta^{p,q}_{m,n}$,
%
\begin{eqnarray}
&& \mathbb{P} \biggl[X_p \geq \frac{\varphi(x)}{m-1}(X_{m+1}+
\cdots+X_n+x),X_q \geq \frac{\varphi(x)}{m-1}(X_{m+1}+
\cdots+X_n+x) \biggr] \nonumber
\\
&& \quad\leq\mathbb{P} \Biggl[\sum_{i=1}^n
w_i\log X_i \geq(w_p+w_q)\log
\frac{\varphi(x)}{m-1} + \log x \Biggr]
\nonumber
\\[-8pt]
\label{estim2prob}
\\[-8pt]
\nonumber
&& \quad\leq N \biggl(-\frac{(w_p+w_q)\log{\varphi(x)}/({m-1}) + \log x -
\sum_{i=1}^n w_i \mu_i}{\sqrt{{\ve{w}}^\perp\cm{\ve{w}}}} \biggr)\notag
\\
&&\quad = \exp \biggl\{-\frac{\log^2 x }{2{\ve{w}}^\perp\cm{\ve{w}}} + \mathrm{O}(\log x \cdot\log\log x) \biggr
\}\nonumber
\end{eqnarray}
as $x\rightarrow\infty$.

Since the matrix $\cm$ is invertible and positive definite, the
mapping ${\ve{w}}\mapsto{\ve{w}}^\perp\cm{\ve{w}}$ is strictly convex.
This implies
that
\[
\min_{{\ve{w}}\in\Delta^{p,q}_{m,n} }{\ve{w}}^\perp\cm{\ve{w}}< \max \Bigl(
\min_{{\ve{w}}\in\Delta^{p}_{m,n} }{\ve{w}}^\perp\cm{\ve{w}},\min
_{{\ve{w}}\in\Delta^{q}_{m,n} }{\ve{w}}^\perp\cm{\ve{w}}\Bigr).
\]
Since
\[
\delta_{4,p} = {\min_{{\ve{w}}\in\Delta^p_{m,n}}{\ve{w}}^\perp
\cm{\ve{w}}},
\]
we conclude from the estimate in \eqref{estim2prob} that the terms in the
second sum in formula \eqref{upperXm} provide a negligible contribution
to the
asymptotics, so that
\begin{eqnarray*}
\mathbb{P}\bigl[X^{(m)}\geq x\bigr] &\leq & \sum
_{p=1}^m\delta_{1,p}(\log
x)^{\delta
_{2,p}}x^{\delta_{3,p}} \exp \biggl\{-\frac{\log^2x}{2 \delta_{4,p}} \biggr\}\bigl(1+\mathrm{O}
\bigl((\log x)^{-1} \bigr)\bigr)
\\
& =&  \delta_{1}(\log x)^{\delta_{2}}x^{\delta_{3}} \exp \biggl
\{-\frac{\log^2 x}{\delta_4} \biggr\}\bigl(1+\mathrm{O} \bigl((\log x)^{-{1}/{2}} \bigr)\bigr)
\end{eqnarray*}
as $x\rightarrow\infty$.

\textit{Lower estimate}.
Let $F_p$, $0\le p\le m$, be random variables. Then for every such $p$,
the following inclusion is valid:
%
\begin{equation}\label{Ekr2}
\{F_p\ge F_0\}\subset\bigcup
_{q\neq p}\{F_p\ge F_0,F_q
\ge F_0\}\cup\{ F_p\ge F_0,F_q<
F_0 \mbox{ for all } q\neq p\}.
\end{equation}
In addition, the sets $\{F_p\ge F_0,F_q< F_0
\mbox{ for all }   q\neq p\}$, $1\le p\le m$, are disjoint. Now, setting
$F_p=X_p$, $1\le p\le m$, and $F_0=X_{m+1}
+ \cdots+ X_n + x$ in (\ref{Ekr2}),\ we can easily derive the
following lower bound for the probability of our interest:
\begin{eqnarray*}
\mathbb{P}\bigl[X^{(m)}\geq x\bigr] &\geq & \sum
_{p=1}^m \mathbb{P}[X_p \geq
X_{m+1} + \cdots+ X_n + x]
\\
&&{}- \sum_{1\leq p,q,\leq m, p\neq q}\mathbb{P}[X_p \geq
X_{m+1} + \cdots+ X_n + x,X_q \geq
X_{m+1} + \cdots+ X_n + x ].
\end{eqnarray*}
Similarly to the first part of the proof, we can now show that the
terms in the second line make a negligible contribution to the limit.
It follows that
\[
\mathbb{P}\bigl[X^{(m)}\geq x\bigr]\geq\delta_{1}(\log
x)^{\delta_{2}}x^{\delta_{3}} \exp \biggl\{-\frac{\log^2 x}{2\delta_4} \biggr
\}\bigl(1+\mathrm{O} \bigl((\log x)^{-{1}/{2}} \bigr)\bigr)
\]
as $x\rightarrow\infty$.
\end{pf*}

\section{Numerics}\label{num.sec}
\subsection{Implementation and domain of validity of the asymptotic formulas}
Formulas \eqref{genw.eq} and \eqref{Egenerass}, as well as formula
\eqref{Eforr1} have a vertical asymptote at $x=1$, which means that
these formulas are not valid for $x\geq1$ and in practice have very
poor accuracy unless $x$ is much smaller than one, which may correspond
to very small probabilities. To partially alleviate this difficulty, we
suggest to use for numerical computations the following natural
modification of formula \eqref{genw.eq}, where the asymptote is shifted
towards the center of the distribution:
%
\begin{eqnarray}
\mathbb{P}[X \leq x] &=& \widetilde C \biggl(\log\frac{1}{x} + \bar{
\ve{w}}^\perp\veee{\mu}+ \mathcal E(\bar{\ve{w}}) \biggr)^{-({1+\bar{n}})/{2}}
\exp \biggl\{-\frac{(\log x-\bar{\ve{w}}^\perp\veee{\mu}
-\mathcal E(\bar{\ve{w}}))^2}{2 \bar{\ve{w}}^\perp\cm\bar{\ve{w}}} \biggr\}
\nonumber
\\[-8pt]
\label{genwshift.eq}
\\[-8pt]
\nonumber
&&{}\times \bigl(1+\mathrm{O} \bigl(|\log x|^{-1} \bigr) \bigr)\qquad  \mbox{as }x\to
0.
\end{eqnarray}
This formula clearly has the same asymptotics as \eqref{genw.eq},
however its domain of validity is larger and numerical experiments show
that it has a better performance. The expression for the distribution
is well defined for all $x$ such that
\[
x \leq x^*= \mathrm{e}^{\bar{\ve{w}}^\perp\veee{\mu}+ \mathcal E(\bar{\ve{w}})}.
\]
For example, assume that $Y_1,\ldots,Y_n$ are identically distributed
with variance $\sigma^2$ and constant correlation $\rho$. Table~\ref{assympt.tab} gives the values of the probability $\mathbb{P}[X\leq
x^*]$ for different values of $\sigma,\rho$ and $n$.
Formulas \eqref{Eforr1}, \eqref{dfdiff.eq} and \eqref{densdiff.eq} can
be modified in a similar manner. Note that the formulas for the
conditional law of Corollaries \ref{laplace.cor}, \ref{laplacedens.cor}
and \ref{laplace2.cor} do not contain a vertical asymptote due to the
cancellation of the logarithmic singularities.

\begin{table}
\tablewidth=\textwidth
\caption{Location of the vertical asymptote in formula \protect\eqref{genwshift.eq} for various values of model parameters}
\label{assympt.tab}
\begin{tabular*}{\tablewidth}{@{\extracolsep{\fill}}lllllllll@{}}
\hline
$\sigma$                & $0.3$  & $1.0$  & $0.3$  & $1.0$  & \phantom{0}$0.3$  & \phantom{0}$1.0$  & \phantom{0}$0.3$  & \phantom{0}$1.0$   \\
$\rho$                  & $0.5$  & $0.5$  & $0.0$  & $0.0$  & \phantom{0}$0.5$  & \phantom{0}$0.5$  & \phantom{0}$0.0$  & \phantom{0}$0.0$   \\
$n$                     & $5$    & $5$    & $5$    & $5$    & $20$   & $20$   & $20$   & $20$    \\
$\mathbb{P}[X\leq x^*]$ & $0.47$ & $0.40$ & $0.40$ & $0.24$ & \phantom{0}$0.46$ & \phantom{0}$0.38$ & \phantom{0}$0.27$ & \phantom{0}$0.041$ \\
\hline
\end{tabular*}
\end{table}

\subsection{Using the asymptotic formulas directly}\label{num1.sec}
\begin{figure}

\includegraphics{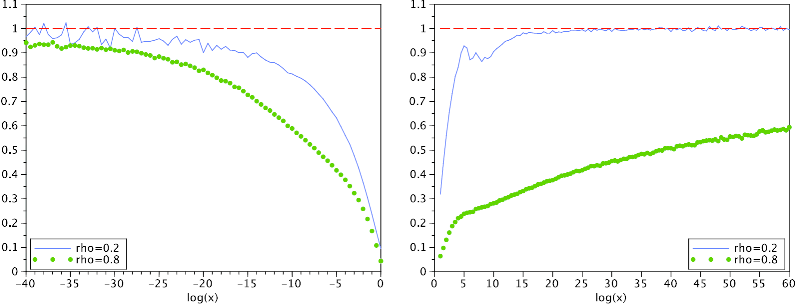}

\caption{Ratios of the Monte Carlo estimate of the distribution
function (survival function) to the estimate obtained using the
asymptotic formulas. Left: $\mathbb{P}[X \leq x]$. Right: $\mathbb{P}[X^{(2)} \geq x]$. The fluctuations in the curves are due to the Monte
Carlo error.}
\label{2cor.fig}
\end{figure}

To
illustrate the performance of the asymptotic formulas of Theorems~\ref{gen.thm} and \ref{Tpresent} numerically, we have taken a $4\times4$
covariance matrix with the following entries: $b_{ij} = \sigma_i \sigma_j \rho$ (constant correlation) where $\sigma= \{2, 2.3,3, 3\}$. The
distribution functions
\[
\mathbb{P}[X \leq x] = \mathbb{P}\bigl[\mathrm{e}^{Y_1} + \mathrm{e}^{Y_2} +
\mathrm{e}^{Y_3} + \mathrm{e}^{Y_4} \leq x\bigr]
\]
and
\[
\mathbb{P}\bigl[X^{(2)} \geq x\bigr] = \mathbb{P}\bigl[\mathrm{e}^{Y_1}
+ \mathrm{e}^{Y_2} - \mathrm{e}^{Y_3} - \mathrm{e}^{Y_4}\geq x\bigr]
\]
have been computed first, using the asymptotic formulas given in
Theorems \ref{gen.thm} and \ref{Tpresent}, with the modification
suggested in formula \eqref{genwshift.eq}. The corresponding asymptotic
approximations will be denoted below by $F_{\mathrm{a}}(x)$ and $F^{(2)}_{\mathrm{a}}(x)$,
respectively. Then we evaluated Monte Carlo estimates $F_{\mathrm{mc}}(x)$ and
$F^{(2)}_{\mathrm{mc}}(x)$ of these quantities (the Monte Carlo algorithm is
described in detail later in this section). To evaluate the quality of
the asymptotic approximation, we plot the ratios $\frac
{F_{\mathrm{mc}}(x)}{F_{\mathrm{a}}(x)}$ and $\frac{F^{(2)}_{\mathrm{mc}}(x)}{F^{(2)}_{\mathrm{a}}(x)}$ for a
wide range of values of $x$. These ratios, plotted as functions of
$\log x$, are shown in Figure~\ref{2cor.fig} for two values of the
correlation coefficient $\rho$.

In the evaluation of the asymptotic formula for $\mathbb{P}[X \leq x]$,
one needs to solve the quadratic programming problem formulated in
\eqref{minprob}. For the first value, $\rho=0.2$, the solution to this
problem is $\bar{\ve{w}} \approx\{0.44\ 0.30\ 0.13\ 0.13\}$. Thus,
here we are in the setting of the ``special case'', where the
asymptotics is obtained directly by Laplace's method (see Lemma~\ref
{easy.prop}). For the second value, $\rho=0.8$, the solution is $\bar
{\ve{w}} \approx\{0.83\ 0.17\ 0\ 0\}$, so only the first two components
make a contribution to the asymptotics.

In the evaluation of the asymptotic formula for $\mathbb{P}[X^{(2)} \geq
x]$, one needs to solve the problem in~\eqref{minprobp} twice, for
$p=1$ and $p=2$, and compare the resulting minimum values. Here, for
$\rho=0.2$, the solutions are $\bar{\ve{w}}_1 = \bar{\ve{w}}_2 = \{1\
0\ 0\}$, and $p=2$ gives a larger minimum value, so that the asymptotic
behavior of the distribution function is determined by the second
component of the vector~$Y$ only. For $\rho= 0.8$, the solutions are
$\bar{\ve{w}}_1 \approx\{1.32\ {-}0.16\ {-}0.16\}$ and $\bar{\ve{w}}_2
\approx\{1.1\ {-}0.05 \ {-}0.05\}$, and once again, the minimum value is
greater for $p=2$. Therefore, in this case the asymptotic behavior is
determined by the second, third and fourth components of $\ve Y$.

Analyzing Figure~\ref{2cor.fig}, one can make the following
observations, which turn out to be rather generic:
\begin{itemize}
\item As expected, the ratio of the distribution functions converges to
one, but this convergence is very slow. This
observation is consistent with the logarithmic error bounds in Theorems
\ref{gen.thm} and \ref{Tpresent}.
\item Although the ratio of the estimates converges to one very slowly,
this ratio is never very far from one (compared to the value of the
probability itself), which means that the asymptotic formula gives the
right order of magnitude for a wide range of probabilities. For
instance, for $\rho=0.8$, the values of $x$, shown in the left graph,
correspond to the range of probabilities from ${\sim}5\times10^{-93}$
for $\log x = -40$ to $0.2$ for $\log x = 0$.
\end{itemize}

\subsection{Efficient Monte Carlo estimation of tail probabilities}
\label{num2.sec}
As we have already seen, due to the slow convergence, the asymptotic
formulas in Theorems \ref{gen.thm} and \ref{Tpresent} typically
provide only order-of-magnitude approximations of the distribution
function of the sum/difference of log-normal random variables. When a
more precise estimate is needed, and the dimension $n$ is large, one
can use a Monte Carlo estimator. In such a case, as we will next
explain, the asymptotic formulas can be utilized to construct very
efficient variance reduction procedures. To save space, we will only
discuss the case of distribution functions. Similar ideas can be used
to reduce the variance of Monte Carlo estimates of densities,
conditional expectations or other quantities of interest.
\subsubsection*{Left tail of $X$}
For the distribution function
$F(x) = \mathbb{P}[X\leq x]$, the standard estimate
is the following:
%
\begin{equation}\label{standard.eq}
\widehat F_N(x) = \frac{1}{N}\sum
_{k=1}^N \mathbf{1}_{\{\sum_{i=1}^n\exp
\{Y_i^{(k)}\}\leq x\}},
\end{equation}
where $\ve Y^{(1)},\ldots,\ve Y^{(N)}$ are i.i.d. vectors with the law
$N(\veee{\mu},\cm)$. However, this estimate is not a suitable
approximation of the tail of the distribution function. Indeed, the
variance of $\widehat F_N(x)$ is given by
\[
\operatorname{Var} \widehat F^{}_N(x) = \frac{F(x)-F^2(x)}{N}
\sim\frac
{F(x)}{N}, \qquad x\to 0,
\]
and the relative error, that is,
\[
\frac{\sqrt{\operatorname{Var}  \widehat F^{}_N(x) }}{F(x)} \sim\frac{1}{\sqrt
{N F(x)}},
\]
explodes very quickly as $x\to0$ (it behaves like $\mathrm{e}^{c\log^2 x}$ for
some constant $c$). The usual way to reduce variance in the Gaussian
context is via importance sampling. The idea is to rewrite the formula
for $F$ as follows:
\[
F(x) = \mathbb E\biggl[\exp\biggl\{-\veee\Lambda^\perp
\cm^{-1}(\ve Y-\veee{\mu}) - \frac
{1}{2} \veee\Lambda^\perp
\cm^{-1}\veee\Lambda\biggr\} \mathbf{1}_{\{\sum_{i=1}^n
\exp\{Y_i + \Lambda_i\}\leq x\}}\biggr],
\]
where $\veee\Lambda\in\mathbb{R}^n$ is a vector that will be chosen
later. Note that if $\veee\Lambda= 0$, then the standard estimate is
recovered. The goal is to find a nonzero $\veee\Lambda$ such that the
corresponding estimate
%
\begin{equation}\label{reduced.eq}
\widehat F^{\veee\Lambda}_N(x) = \frac{1}{N}\sum
_{k=1}^N \exp\biggl\{-{\veee\Lambda}^\perp
\cm^{-1}\bigl(Y^{(k)}- \veee{\mu} \bigr) - \frac{1}{2} {
\veee\Lambda }^\perp\cm^{-1}{\veee\Lambda} \biggr\}
\mathbf{1}_{\{\sum_{i=1}^n\exp\{
Y^{(k)}_i + \Lambda_i\}\leq x\}}
\end{equation}
has variance smaller than that of the standard estimate.

Simple computations show that the variance of $\widehat F^{{\veee\Lambda
}}_N(x)$ is given by
\begin{eqnarray*}
\operatorname{Var} \widehat F^{{\veee\Lambda}}_N(x) &=& \frac{1}{N}
\operatorname{Var} \biggl[\exp\biggl\{-{\veee\Lambda}^\perp\cm^{-1}({\ve
Y}-\veee{\mu}) - \frac{1}{2} {\veee\Lambda}^\perp
\cm^{-1}{\veee\Lambda} \biggr\} \mathbf{1}_{\{\sum_{i=1}^n
\exp\{Y_i + \Lambda_i\}\leq x\}}\biggr]
\\
&=& \frac{1}{N} \bigl\{\mathbb E\bigl[\exp\bigl\{-2{\veee\Lambda}^\perp\cm^{-1}({\ve Y}-\veee{\mu}) - {\veee\Lambda
}^\perp\cm^{-1}{\veee\Lambda} \bigr\} \mathbf{1}_{\{\sum_{i=1}^n \exp\{Y_i + \Lambda_i\}\leq x\}}
\bigr] - F^2(x) \bigr\}
\\
& =& \frac{1}{N} \biggl\{\exp\biggl\{\frac{1}{2} {\veee\Lambda}^\perp\cm ^{-1}{\veee\Lambda}\biggr\} \mathbb E\bigl[
\exp\bigl\{-{\veee\Lambda}^\perp\cm^{-1}({\ve Y}-\veee{\mu})
\bigr\} \mathbf{1}_{\{\sum_{i=1}^n \exp\{Y_i\}\leq x\}}\bigr] - F^2(x) \biggr\}
\\
& =& \frac{1}{N} \Biggl\{\exp\bigl\{{\veee\Lambda}^\perp
\cm^{-1}{\veee\Lambda}\bigr\} \mathbb{P} \Biggl[\sum
_{i=1}^n\exp\{Y_i-\Lambda_i
\}\leq x \Biggr] - F^2(x) \Biggr\}.
\end{eqnarray*}
Let
\[
V({\veee\Lambda},x) =\exp\bigl\{{\veee\Lambda}^\perp\cm^{-1}{\veee\Lambda}\bigr\} \mathbb{P} \Biggl[\sum_{i=1}^n
\mathrm{e}^{Y_i-\Lambda_i }\leq x \Biggr].
\]
Since $F(x)$ does not depend on ${\veee\Lambda}$, the optimal variance
reduction is obtained by minimizing $V({\veee\Lambda},x)$ as a function
of ${\veee\Lambda}$. Our idea is to obtain an explicit estimate by
replacing the probability in the previous expression by an
asymptotically equivalent expression given in Theorem~\ref{gen.thm}. In
other words, we compute an approximation to the optimal ${\veee\Lambda}$
by minimizing
\begin{eqnarray*}
\widetilde V({\veee\Lambda}, x) &=& {\veee\Lambda}^\perp\cm^{-1}{
\veee\Lambda} -\frac{1}{2}\sum_{i,j=1}^{\bar{n}}
\bar a_{ij} \biggl(\log\frac
{\bar A_1+\cdots+ \bar A_{\bar{n}}}{\bar A_i} + \bar\mu_i -
\log x - \Lambda_i \biggr)
\\
&&{}\times \biggl(\log\frac{\bar A_1+\cdots+ \bar
A_{\bar{n}}}{\bar A_j} + \bar\mu_j -\log x -
\Lambda_j \biggr).
\end{eqnarray*}
To obtain the expression above, we have omitted all the factors in the
formula in Theorem~\ref{gen.thm}, which do not depend on ${\veee\Lambda
}$, and have also taken the logarithm of the resulting expression. The
optimal value ${\veee\Lambda}^*$ of ${\veee\Lambda}$ can be found by
solving the following system of equations:
\[
\frac{\partial\widetilde{V}}{\partial\Lambda_i}\bigl({\veee\Lambda }^*\bigr)=0, \qquad 1\le i\le n.
\]
This system can be rewritten as follows:
\begin{eqnarray*}
2\sum_{j=1}^na_{ij}
\Lambda^*_j+\sum_{j=1}^{\bar{n}}
\bar a_{ij} \biggl(\log\frac{\bar A_1+\cdots+ \bar A_{\bar{n}}}{
\bar A_j} + \bar\mu_j -
\log x - \Lambda_j \biggr) &=&0, \qquad  1\le i\le \bar{n},
\\
2\sum_{j=1}^na_{ij}
\Lambda^*_j &=& 0,\qquad  \bar{n}< i\le n.
\end{eqnarray*}
Applying the matrix $\cm$ to the previous system, we obtain
%
\begin{equation} \label{Eadd}
2\Lambda^*_k + \sum_{i,j=1}^{\bar{n}}
b_{ki} \bar a_{ij} \biggl(\log\frac
{\bar A_1+\cdots+ \bar A_{\bar{n}}}{\bar A_j} + \bar
\mu_j -\log x - \Lambda^*_j \biggr)=0,
\end{equation}
for all $k=1,\ldots,n$. When $k\leq\bar{n}$, the formula in (\ref{Eadd}) simplifies to
\[
\Lambda^*_k + \log\frac{\bar A_1+\cdots+ \bar A_{\bar{n}}}{\bar A_k} + \bar\mu_k -
\log x=0.
\]
Substituting this into (\ref{Eadd}), we see that for all $k$, the
optimal value $\Lambda^*_k$ is given by
%
\begin{equation}\label{lambdaopt.eq}
\Lambda^*_k = \sum_{i,j=1}^{\bar{n}}
b_{ki} \bar a_{ij} \biggl(\log x - \log\frac{\bar A_1+\cdots+ \bar A_{\bar{n}}}{\bar A_j}
- \bar\mu_j \biggr).
\end{equation}
Note that since the optimal vector ${\veee\Lambda}^*$ depends on $x$, we
cannot apply Theorem~\ref{gen.thm} directly to characterize the
asymptotic behavior of the function $V({\veee\Lambda}^*,x)$ as $x\to0$.
Nevertheless,\vadjust{\goodbreak} this function can be estimated from above by using
Jensen's inequality as follows:
\begin{eqnarray*}
V\bigl({\veee\Lambda}^*,x\bigr) &\leq & \mathrm{e}^{{\veee\Lambda}^{*\perp} \cm^{-1}{\veee\Lambda}^*}\mathbb{P} \Biggl[\sum
_{i=1}^n \bar{w}_i
\bigl({Y_i-\Lambda_i^* - \log\bar{w}_i}\bigr)
\leq\log x \Biggr]
\\
&=& \mathrm{e}^{{\veee\Lambda}^{*\perp} \cm^{-1}{\veee\Lambda}^*} N \biggl(\frac{\bar
{\ve{w}}^\perp{\veee\Lambda}^* + \log x + \sum_{i=1}^n \bar{w}_i(\log
\bar{w}_i - \mu_i)}{\sqrt{\bar{\ve{w}}^\perp\cm\bar{\ve{w}}}} \biggr),
\end{eqnarray*}
where $\bar{\ve{w}}$ is the solution of \eqref{minprob} and $N$ is the
standard normal distribution function. Substituting the expression in
\eqref{lambdaopt.eq} for ${\veee\Lambda}^*$ and using \eqref{awstar}, we obtain
\[
\bar{\ve{w}}^\perp{\veee\Lambda}^* = \log x + \sum
_{i=1}^n \bar{w}_i(\log\bar{w}_i - \mu_i)
\]
and
\begin{eqnarray*}
{\veee\Lambda}^{*\perp} \cm^{-1}{\veee\Lambda}^* &=&   \sum
_{i,j=1}^{\bar{n}} \bar a_{ij} (\log x +\log\bar{w}_i - \bar\mu_i ) (\log x +\log\bar{w}_j -
\bar\mu_j )
\\
& =&  \frac{1}{\bar{\ve{w}}^\perp\cm\bar{\ve{w}}} \Biggl\{\log x + \sum_{i=1}^n
\bar{w}_i(\log\bar{w}_i - \mu_i) \Biggr
\}^2
\\
&&{}+ \sum_{i,j=1}^{\bar{n}} \bar a_{ij}
(\log\bar{w}_i - \bar\mu_i ) (\log\bar{w}_j -
\bar\mu_j ) - \frac{1}{\bar{\ve{w}}^\perp
\cm\bar{\ve{w}}} \Biggl\{\sum
_{i=1}^n \bar{w}_i(\log\bar{w}_i - \mu _i) \Biggr\}^2.
\end{eqnarray*}
Therefore, as $x\to0$,
\[
V\bigl({\veee\Lambda}^*,x\bigr) \lesssim C\frac{\exp \{-({1}/({\bar{\ve{w}}^\perp\cm\bar{\ve{w}}})) \{\log x + \sum_{i=1}^n \bar{w}_i(\log
\bar{w}_i - \mu_i) \}^2 \}}{\log{1}/{x}},
\]
where the constant $C$ is independent of $x$. Comparing the previous
estimate with the asymptotics of $F(x)$ (see formula \eqref{genw.eq}),
we see that, for a different constant $C$,
%
\begin{equation}
V\bigl({\veee\Lambda}^*,x\bigr) \lesssim C F^2(x) \biggl(\log
\frac{1}{x} \biggr)^{\bar{n}}
\end{equation}
as $x\to0$. This means in particular that our estimator is \emph
{logarithmically efficient} in the sense of Asmussen and Glynn \cite{AG}, Section VI.1.

To test the performance of the proposed variance reduction algorithm,
we have computed the Monte Carlo estimates with and without variance
reduction for different levels $x$, using the same numerical values of
the parameters as above. Table~\ref{ratio.tab} shows the relative error
of the estimate~\eqref{reduced.eq}, where the value of ${\veee\Lambda
}^*$ is given by \eqref{lambdaopt.eq} as well as the ratio of the \emph
{standard deviation} of the estimate \eqref{standard.eq} to that of the
estimate \eqref{reduced.eq}. The relative errors appear quite stable
and the reduction factors are greater than one for all values of $x$
and in general quite spectacular, ranging from 4--5 for not so small
probabilities of order of $1\%$ to hundreds for probabilities of order
of $10^{-6}$.

\begin{table}
\tabcolsep=0pt
\tablewidth=\textwidth
\caption{Relative errors of the variance reduction estimate \protect\eqref{reduced.eq} with the value of ${\veee\Lambda}^*$ given by
\protect\eqref{lambdaopt.eq}, and the factors by which the \emph{standard deviation}
of the estimate is reduced with the variance reduction algorithm}
\label{ratio.tab}
\begin{tabular*}{\tablewidth}{@{\extracolsep{\fill}}llllllll@{}}
\hline
\multicolumn{4}{@{}l}{$\rho= 0.2$} & \multicolumn{4}{l@{}}{$\rho=0.8$} \\[-4pt]
\multicolumn{4}{@{}l}{\hrulefill}    & \multicolumn{4}{l@{}}{\hrulefill}\\
$x$ & $\mathbb{P}[X\leq x]$ & rel. error& red. factor & $x$ & $\mathbb{P}[X\leq x]$ & rel. error & red. factor\\
\hline
0.006738 & 0.0000027 & 0.43\% & 152.8 & 0.0002035 & 0.0000012 & 0.27\%
&269 \\
0.01831& 0.0000424 & 0.37\% & \phantom{0}38.07 & 0.0009119 & 0.0000331 & 0.26\% &
\phantom{0}69.08 \\
0.04979 & 0.0004639 & 0.32\%& \phantom{0}14.48 & 0.004089 & 0.0005282 &0.25\% &
\phantom{0}16.07 \\
0.1353& 0.003457 & 0.28\% & \phantom{00}6.188 & 0.01832 & 0.005085 & 0.25\% &\phantom{00}5.312
\\
0.3679 & 0.01798 & 0.24\%& \phantom{00}3.152& 0.08209 & 0.02998 & 0.26\%& \phantom{00}2.256 \\
1 & 0.06603 & 0.20\% &\phantom{00}1.845 & 0.3679 & 0.1141 & 0.27\%& \phantom{00}1.078 \\
\hline
\end{tabular*}
\end{table}

\begin{figure}

\includegraphics{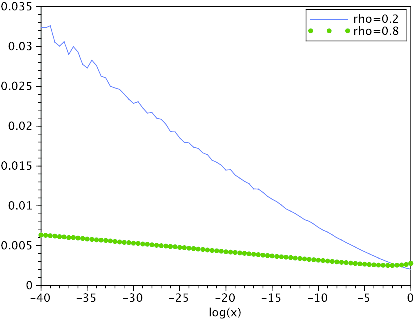}

\caption{Relative error of the variance reduction estimate \protect\eqref
{reduced.eq} with the value of ${\veee\Lambda}^*$ given by \protect\eqref{lambdaopt.eq}.}
\label{relerr.fig}
\end{figure}

Figure~\ref{relerr.fig} shows the relative error of the estimate \eqref
{reduced.eq} with the value of ${\veee\Lambda}^*$ given by \eqref
{lambdaopt.eq} (the standard deviation divided by the value of the
estimate, computed over $10^6$ trajectories). As shown by the
theoretical analysis of the variance, the relative error\vspace*{1pt} grows only
logarithmically in~$x$, which means that even for very small
probabilities (such as $10^{-100}$), our estimator requires a
reasonable number of trajectories to obtain adequate precision.
%

\subsubsection*{Additional tests}
To evaluate the robustness of our variance reduction method with
respect to the choice of the parameters of the model and the number of
variables $n$, we have performed additional tests, assuming this time
that $Y_1,\ldots,Y_n$ are identically distributed with law $N(0,\sigma
^2)$ and have a constant correlation $\rho$. Table~\ref{addtest.tab}
shows the \emph{standard deviation} reduction factors for various
values of $\sigma$, $\rho$ and $n$. For each test, the value of $x$ was
selected so that the probability $\mathbb{P}[X\leq x]$ is approximately
equal to $10^{-3}$ (belongs to the interval $(0.9\times10^{-3},
1.1\times10^{-3})$). We see that in all tests but one, the standard
deviation is reduced by a factor greater than $10$, which means that,
for equal precision, the computation would be accelerated by a factor
greater than $100$.

\begin{table}[b]
\caption{\emph{Standard deviation} reduction factors for the additional
tests. The probability $\mathbb{P}[X\leq x]$ is approximately equal to
$10^{-3}$ for all tests}
\label{addtest.tab}
\begin{tabular*}{\tablewidth}{@{\extracolsep{\fill}}lllllllll@{}}
\hline
$\sigma$ & \phantom{0}$0.3$& \phantom{0}$1.0$ & \phantom{0}$0.3$ & \phantom{0}$1.0$& \phantom{0}$0.3$& \phantom{0}$1.0$ & \phantom{0}$0.3$ & \phantom{0}$1.0$\\
$\rho$ & \phantom{0}$0.5$& \phantom{0}$0.5$& \phantom{0}$0.0$ & \phantom{0}$0.0$& \phantom{0}$0.5$& \phantom{0}$0.5$& \phantom{0}$0.0$ & \phantom{0}$0.0$\\
$n$ & \phantom{0}$5$& \phantom{0}$5$& \phantom{0}$5$ & \phantom{0}$5$& $20$& $20$& $20$ & $20$\\
$x$ & \phantom{0}$2.5$& \phantom{0}$0.55$& \phantom{0}$3.4$ & \phantom{0}$1.6$ & $10.5$ & \phantom{0}$2.7$ & $16.9$ & $14.3$ \\
red.~factor & $15.7$ & $14.7$ & $14.0$ & $10.1$ & $15.2$ & $14.2$ &
$11.4$ & \phantom{0}$4.8$ \\
\hline
\end{tabular*}
\end{table}

\subsubsection*{Right tail of a log-normal difference}
In this case, the standard estimate of the survival function has the form
%
\begin{equation}\label{plaindif.eq}
\widehat F_N(x) = \frac{1}{N}\sum
_{k=1}^N \mathbf{1}_{\{\sum_{i=1}^m \exp
\{Y_i^{(k)}\} -
\sum_{i=m+1}^n \exp\{Y_i^{(k)}\} \geq x\}},
\end{equation}
and the alternative estimate which may potentially reduce variance is
as follows:
%
\begin{eqnarray}
\widehat F^{\veee\Lambda}_N(x) &=&  \frac{1}{N}\sum
_{k=1}^N\exp\biggl\{-{\veee\Lambda}^\perp
\cm^{-1}\bigl({\ve Y}^{(k)}-\veee{\mu}\bigr) -
\frac{1}{2} {\veee\Lambda}^\perp\cm^{-1}{\veee\Lambda}
\biggr\}
\nonumber
\\[-8pt]
\label{reduced2.eq}
\\[-8pt]
\nonumber
&&{}\times\mathbf{1}_{\{\sum_{i=1}^m\exp\{Y_i^{(k)} + \Lambda_i\} -
\sum_{i=m+1}^n\exp\{Y_i^{(k)}+\Lambda_i\}
\geq x\}}.
\end{eqnarray}
To find the optimal value of ${\veee\Lambda}$, we need to minimize
\[
\exp\bigl\{{\veee\Lambda}^\perp\cm^{-1}{\veee\Lambda}\bigr\}
\mathbb{P} \Biggl[\sum_{i=1}^m\exp\bigl
\{Y_i^{(k)} - \Lambda_i\bigr\} - \sum
_{i=m+1}^n \exp\bigl\{Y_i^{(k)}-
\Lambda_i\bigr\} \geq x \Biggr],
\]
and once again, the main idea is to minimize the asymptotic
approximation to this function, given in Theorem~\ref{Tpresent}.
Assuming for simplicity that the set $\mathcal P_4$ defined in \eqref{P4} is a singleton, $\mathcal P_4 = \{p\}$, the problem reduces to
that of minimizing the following function:
\begin{eqnarray*}
\widetilde V({\veee\Lambda},x) &=&  {\veee\Lambda}^\perp\cm^{-1} {
\veee\Lambda} -\frac{1}{2}\sum_{i,j=1}^{\bar{n}^{(p)}}
\bar a^{(p)}_{ij} \biggl(\log\frac{\bar A^{(p)}_1+\cdots+ \bar A^{(p)}_{\bar{n}^{(p)}}}{|\bar A^{(p)}_i|} + \bar
\mu^{(p)}_i -\log x - \Lambda_{\bar
k^{(p)}(i)} \biggr)
\\
&&{}\times \biggl(\log\frac{\bar A^{(p)}_1+\cdots+
\bar A^{(p)}_{\bar{n}^{(p)}}}{|\bar A^{(p)}_j|} + \bar\mu^{(p)}_j -
\log x - \Lambda_{\bar k^{(p)}(j)} \biggr).
\end{eqnarray*}
Next, reasoning as in the proof of (\ref{lambdaopt.eq}), we see that
the optimal value ${\veee\Lambda}^*$ of ${\veee\Lambda}$ is given by
%
\begin{equation}\label{lambdaopt2.eq}
\Lambda^*_k = \sum_{i,j=1}^{\bar{n}^{(p)}}
b_{k,\bar k^{(p)}(i)} \bar a^{(p)}_{ij} \biggl(\log x - \log
\frac{\bar A^{(p)}_1+\cdots+ \bar
A^{(p)}_{\bar{n}^{(p)}}}{|\bar A^{(p)}_j|} - \bar\mu^{(p)}_j \biggr).
\end{equation}
However, here the computation remains only heuristic, since there is no
simple upper bound for the variance of the estimator with the optimal
${\veee\Lambda}^*$.

\begin{table}[b]
\caption{Relative errors of the variance reduction estimate \protect\eqref{reduced2.eq}
with the value of ${\veee\Lambda}^*$ given by \protect\eqref{lambdaopt2.eq}, and those of the plain Monte Carlo estimate
\protect\eqref{plaindif.eq}, whenever available}
\label{ratio2.tab}
\begin{tabular*}{\tablewidth}{@{\extracolsep{\fill}}llllllll@{}}
\hline
\multicolumn{4}{@{}l}{$\rho= 0.2$} & \multicolumn{4}{l@{}}{$\rho=0.8$} \\[-4pt]
\multicolumn{4}{@{}l}{\hrulefill} &  \multicolumn{4}{l@{}}{\hrulefill}\\
$x$ & $\mathbb{P}[X^{(2)}\geq x]$ & rel. error&  rel. error,   & $x$ & $\mathbb{P}[X^{(2)}\geq x]$ & rel.
error & rel. error,\\
&&&plain MC    &&&&    plain MC\\
\hline
$\mathrm{e}$ & 0.2672 & \phantom{0}0.153\%& \phantom{0}0.166\% & $\mathrm{e}$ & 0.1121 &0.74\%& \phantom{00}0.281\% \\
$\mathrm{e}^5$ & 0.01564 &\phantom{0}2.37\% & \textup{0}0.792\% & $\mathrm{e}^5$ & $2.134\times10^{-3}$
&4.26\%& \phantom{00}2.17\% \\
$\mathrm{e}^{10}$ & $6.771\times10^{-6}$ & 57.1\%& 29.5\% & $\mathrm{e}^{10}$ &
$3.759\times10^{-7}$ &1.15\%& 376\% \\
$\mathrm{e}^{15}$ & $3.459\times10^{-11}$ & \phantom{0}0.274\%& -- & $\mathrm{e}^{15}$ &
$9.765\times10^{-13}$ &0.44\% & -- \\
$\mathrm{e}^{20}$ & $1.724\times10^{-18}$ & \phantom{0}0.318\%& -- & $\mathrm{e}^{20}$ &
$2.654\times10^{-20}$ & 0.502\%& -- \\
$\mathrm{e}^{25}$ & $8.050\times10^{-28}$ & \phantom{0}0.358\% & -- & $\mathrm{e}^{25}$ &
$6.872\times10^{-30}$ &0.561\%& -- \\ \hline
\end{tabular*}
\end{table}

Numerical tests (see Table~\ref{ratio2.tab}) are much less conclusive
than those for the left tail presented in the previous paragraph. For
moderate values of $x$, the algorithm does not lead to any variance
reduction and may even increase variance. However, very far in the
tail, when the probability in question is so small that it cannot be
computed with the conventional Monte Carlo estimator in reasonable
time, the variance reduction estimator becomes very efficient. We
conclude that for the case of log-normal differences, our variance
reduction algorithm can be potentially very useful for the simulation
of extremely rare events (with probability smaller than $10^{-6}$), but
further research and further improvements to the algorithm are
necessary before it can be used in the context of not-so-rare events,
such as those with probability of $10^{-2}$--$10^{-3}$ arising, for
example, in the Value at Risk calculations in financial risk management.

\section{Risk management in the multidimensional Black--Scholes model}
\label{rm.sec}

The tail estimates obtained in this paper can be applied to risk management
problems in the context of the $n$-dimensional Black--Scholes model.
Suppose that the dynamics
of the asset price vector $\ve S_t = (S^1,\ldots,S^n)$ is described by
the following $n$-dimensional stochastic process:
%
\begin{equation}\label{EBSS}
\log S_t=\log\ve S_0+\veee\theta t-
\frac{ \operatorname{diag}(\cm)t}{2}+\cm ^{{1}/{2}}\ve W_t,
\end{equation}
where $\ve W$ is an $n$-dimensional standard Brownian motion, $\cm$ is
the covariance matrix,
$\veee\theta$ is the drift vector, and $ \operatorname{diag}(\cm)$ stands for
the main diagonal of $\cm$.

Consider a portfolio containing the assets $S^i$, $1\le i\le n$ with
weights $\xi_1,\ldots,\xi_n$,
and the price process $X$ defined by
%
\begin{equation}\label{EBS}
X_t = \sum_{i=1}^n
\xi_i S^i_t,\qquad  t\ge0.
\end{equation}
The initial condition for the process $X$ is given by $X_0= \sum_{i=1}^n \xi_i S^i_0$.
The process $X$ can be alternatively expressed as
$X_t=\sum_{i=1}^n\operatorname{sgn}  \xi_i \exp\{Y_t^i\}$,
where
%
\begin{equation}\label{Epr1}
Y_t^i=\log S_0^i+\log|
\xi_i|+\theta_it-\frac{b_{ii}t}{2}+\sum
_{j=1}^n\gamma_{ij}W_t^j,
\qquad 1\le i\le n.
\end{equation}
In (\ref{Epr1}), the symbols $\gamma_{ij}$ stand for the elements of
the matrix $\cm^{{1}/{2}}$. We also set
%
\begin{equation}\label{Emuq}
\mu_{i,t}=\log S_0^i+\log|
\xi_i|+\theta_it-\frac{b_{ii}t}{2}, \qquad 1\le i\le n.
\end{equation}
%
In the sequel, $t$ will be fixed, and the asymptotic formulas obtained in
Sections~\ref{Snachalo} and \ref{Stail1} will be applied to the
random variable $X_t$ defined in
(\ref{EBS}).
The Gaussian data associated with the case described above are given by
the following: the mean vector
is $\vee{\mu}=(\mu_{1,t},\ldots,\mu_{n,t})$ and the covariance matrix is
$t\cm$.



For the purposes of risk management, it is important to solve two
classes of problems in relation with the portfolio $X$:
\begin{itemize}
\item Quantify the behavior of one portfolio in specific adverse
scenarios of market evolution, which are typically defined in terms of
another portfolio (the benchmark). This can be done using our
characterization of the asymptotic behavior of a Gaussian vector
conditionally on the sum or difference of exponentials of its
components (Corollaries \ref{laplacedens.cor} and \ref{laplace2.cor}).
We address this issue in detail in the following paragraph.
\item Evaluate various risk measures for the portfolio $X$, such as the
probability of loss of a given magnitude or the Value at Risk (the
quantile function). The probability of loss may be approximated using
the asymptotic formulas of Section~\ref{Snachalo} (for portfolios with
only positive weights) or Section~\ref{Stail1} (for portfolios with
both positive and negative weights). The asymptotic behavior of the
Value at Risk when the confidence level tends to one
is characterized in Section~\ref{Siv}.
\end{itemize}

\subsection{Behavior of log-normal portfolios under adverse
scenarios}\label{SSlnpb}
Suppose that an investor holds a portfolio containing assets $S^1,\ldots,
S^n$ with weights $v_1,\ldots,v_n$.
The value of such a portfolio is given by
%
\begin{equation}\label{Eporto}
V_t = \sum_{i=1}^n
v_i S^i_t.
\end{equation}
The 1996 Market Risk Amendment to Basel I \cite{mra}
as well as Basel II and Basel III Capital Accords require banks and
investment firms to conduct stress tests to determine their ability to
respond to adverse market events. These adverse scenarios are typically defined
in terms of the performance of a certain benchmark and correspond to a
stylized version of certain crisis events observed in the past. We will
next describe
some examples of plausible stress scenarios and explain how
the corresponding benchmark process $X$ can be defined.
\begin{itemize}
\item \textit{Equity market fall of a certain magnitude}.  This is the
most common stress scenario. The benchmark process $\{X\}_{t\ge0}$
under such a scenario is
the normalized market index, having the initial value $1$, and the
adverse event is $\{X_t = x\}$ for some $t>0$ and $x$ which is supposed
small. The weights $\xi_i$ are then positive
and equal to the normalized market capitalizations of the stocks.
\item  \textit{A certain difference in performance between the equity markets
of two geographical areas or two sectors}.  For instance, one may
assume that the American markets outperform the European
ones, or that small capitalization shares outperform large
capitalization ones. Let
$X^a_t = \sum_{i=1}^m \xi_i S^i_t$
be the market index of the first area, where $\xi_1,\ldots,\xi_m$ are
the positive market capitalization weights of the stocks
$S^1,\ldots,S^m$, and
$X^b_t = \sum_{i=m+1}^n \xi_i S^i_t$
be the market index of the second area. The stress scenario is the
event
\[
\biggl\{\frac{X^a_t}{X^a_0} - \frac{X^b_t}{X^b_0} = x \biggr\}.
\]
This can be dealt with in our framework by taking
\[
X_t = \sum_{i=1}^m
\frac{\xi_i}{X^a_0} S^i_t - \sum
_{i=m+1}^n \frac{\xi_i}{X^b_0} S^i_t
\]
with the stress scenario $\{X_t = x\}$. Here the value of $x$ is large.
\item \textit{A certain difference in performance between two benchmarks}.
The
investor may be interested, for example, in the event when her
portfolio severely underperforms the market. This is similar to the
case considered
above, except that the two benchmarks may contain the same
stocks. Let the two benchmarks be given by
$X^a_t = \sum_{i=1}^n \xi^a_i S^i_t$ and $X^b_t = \sum_{i=1}^n \xi^b_i S^i_t$.
We are once again interested in the stress scenario
$ \{\frac{X^a_t}{X^a_0} - \frac{X^b_t}{X^b_0} =
x \}$. This is equivalent to taking
\[
X_t = \sum_{i=1}^n \biggl\{\frac{\xi^a_i}{X^a_0} - \frac{\xi
^b_i}{X^b_0} \biggr\} S^i_t
\]
and using the stress scenario $\{X_t = x\}$ with $x$ large.
\end{itemize}

Our next goal is to characterize the asymptotic behavior of various
conditional expected values for the portfolio with the price process
given by (\ref{Eporto}) under the stress scenarios described above.
This can be done for the individual stocks or for the entire portfolio.
In the former case, we approximate the conditional probabilities of the form
%
\begin{equation}\label{Econex2}
e_i(t,x) = \mathbb E\Biggl[S^i_t\Big| \sum
_{k=1}^n \xi_k
S^k_t = x\Biggr],
\end{equation}
while in the latter case we deal with the following conditional probabilities:
%
\begin{equation}\label{Econex1}
\mathbb E[V_t | X_t = x] = \sum
_{i=1}^n v_i e_i(t,x).
\end{equation}
%

The quantities $e_i(t,S)$ can be estimated using formulas \eqref{ELT} and
\eqref{lt2.eq} for the conditional Laplace transform, since
%
\begin{equation}\label{Elap1}
e_i(t,x)=\frac{1}{\xi_i}\mathbb{E}\Biggl[\exp\bigl
\{Y_t^i\bigr\}\Big|\sum_{k=1}^n
\exp\bigl\{ Y_t^k\bigr\}= x\Biggr]
\end{equation}
for all $1\le i\le n$.
The following results are thus direct consequences of Corollaries
\ref{laplace.cor} and \ref{laplace2.cor}.

\begin{theorem}\label{stress.prop}
Suppose that the weights $\xi_1,\ldots,\xi_n$ are positive and that
assumption $(\mathcal A)$ holds for the covariance matrix $\cm$. Then
the following are true:
\begin{enumerate}
\item If $1\le i\le\bar{n}$, then as $x\rightarrow0$,
\begin{eqnarray*}
e_i(t,x) &=&\frac{x}{\xi_i} \frac{\bar A_i}{\bar A_1+\cdots+ \bar A_{\bar{n}}} \biggl(1+\mathrm{O} \biggl(
\biggl(\log \frac{1}{x} \biggr)^{-1} \biggr) \biggr)
\\
& =& x\frac{ \bar{w}_i}{\xi_i} \biggl(1+\mathrm{O} \biggl( \biggl(\log \frac{1}{x}
\biggr)^{-1} \biggr) \biggr).
\end{eqnarray*}
\item If $\bar{n}< i\le n$, then as $x\rightarrow0$,
\begin{eqnarray*}
e_i(t,x)
&=& x^{\sum_{j=1}^{\bar{n}}\bar{A}_jb_{ij}}S_0^i \exp \Biggl\{
\theta_i t- \sum_{p,q=1}^{\bar{n}}
b_{pi} \bar a_{pq} \biggl(\log\frac{\bar A_1+\cdots+ \bar A_{\bar{n}}}{\bar A_q} +
\mu_{q,t} \biggr) \Biggr\}
\\
&&{}\times\exp \Biggl\{-\frac{t}{2}\sum_{p,q=1}^{\bar{n}}
\bar a_{pq} b_{pi} b_{qi} \Biggr\} \biggl(1+\mathrm{O}
\biggl( \biggl(\log \frac{1}{x} \biggr)^{-1} \biggr) \biggr).
\end{eqnarray*}
\end{enumerate}
\end{theorem}

\begin{remark}
Since $\sum_{j=1}^{\bar{n}}\bar{A}_jb_{ij}>1$ for $i\notin\bar I$ (see Remark~\ref{remgt1}),
it follows from Theorem~\ref{stress.prop} that the assets in the
market index can be classified into two
categories, depending on the behavior of their conditional expectation
under the conditional law.
\begin{itemize}
\item``Safe assets'', whose conditional expectations decay
proportionally to the value $x$ of the market
index. Those are exactly the assets, which enter the Markowitz
minimal variance portfolio (solution of problem \eqref{minprob}) with
strictly positive weights.
\item``Dangerous assets'', whose conditional expectations decay faster
than the index.
\end{itemize}
\end{remark}
%

The next assertion concerns the second and the third typical stress
scenarios described above.

\begin{theorem}
Suppose that for $m\leq n$ the weights $\xi_1,\ldots,\xi_m$ are
positive and $\xi_{m+1},\ldots,\xi_n$ are negative, that assumption
$(\mathcal A^i_1)$ holds for matrix $\cm$ with every $i=1,\ldots,m$,
and that the set $\mathcal P_4$ defined in \eqref{P4} is a singleton,
$\mathcal P_4 = \{p\}$. Then the
following are true.
\begin{enumerate}
\item If $i\in I^{(p)}$, then as $x\rightarrow+\infty$,
\begin{eqnarray*}
e_i(t,x) &=& \frac{x}{\xi_i} \frac{|\bar A^{(p)}_i|}{\sum_{j=1}^{\bar{n}^{(p)}} \bar A^{(p)}_j} \times \biggl(1+\mathrm{O}
\biggl( \biggl(\log \frac{1}{x} \biggr)^{-1} \biggr) \biggr)
\\
& =&  x\frac{ |\bar{w}_i|}{\xi_i} \times \biggl(1+\mathrm{O} \biggl( \biggl(\log
\frac{1}{x} \biggr)^{-1} \biggr) \biggr).
\end{eqnarray*}
\item If $i\notin I^{(p)}$, then as $x\rightarrow+\infty$,
\begin{eqnarray*}
e_i(t,x)
&=& S_{0}^ix^{\sum_{j=1}^{\bar{n}^{(p)}}\bar A_j^{(p)}b_{\bar
k^{(p)}(j),i} }\notag
\\
&&{}\times\exp \Biggl\{ \theta_i t-\sum_{j,k=1}^{\bar{n}^{(p)}}
\bar a^{(p)}_{jk} b_{\bar k^{(p)}(j),i} \biggl(\log
\frac{\sum_{l=1}^{\bar{n}^{(p)}} \bar A^{(p)}_l}{|\bar A^{(p)}_k|}+ \mu_{\bar k^{(p)}(k),t} \biggr) \Biggr\} \notag
\\
&&{}\times\exp \Biggl\{ -\frac{t}{2}\sum_{j,k=1}^{\bar{n}^{(p)}}
\bar a^{(p)}_{jk} b_{\bar
k^{(p)}(j),i} b_{\bar k^{(p)}(k),i} \Biggr
\} \biggl(1+\mathrm{O} \biggl( \biggl(\log \frac{1}{x} \biggr)^{-1}
\biggr) \biggr).
\end{eqnarray*}
\end{enumerate}
\end{theorem}

\begin{remark}
It follows from assumption $(\mathcal A^p_1)$ and Remark~\ref{rema1} that for $i\notin I^{(p)}$,
\[
\sum_{j=1}^{\bar{n}^{(p)}}\bar A_j^{(p)}b_{\bar k^{(p)}(j),i}<1.
\]
Therefore, the
assets in the benchmark can once again be classified into the following two
categories:
\begin{itemize}
\item Those assets, whose conditional expectations, given the stress scenario,
grow proportionally to $x$. This category includes exactly one
asset among $S_1,\ldots,S_m$, that one with the highest relative
asymptotic variance with respect to $S_{m+1},\ldots,S_n$. It may or
may not include some assets among $S_{m+1},\ldots,S_n$.
\item Those assets, whose conditional expectations, given the stress scenario,
grow slower than $x$.
\end{itemize}
In other words, the fact that the portfolio $S^1+\cdots+S^m$ strongly outperforms
the portfolio $S^{m+1}+\cdots+S^n$ can be attributed asymptotically to a
very strong performance of a single stock among $S^1,\ldots,S^m$, which
may be partially offset by the performance of some stocks from the
second group.
\end{remark}


\subsection{Log-normal portfolios and Value at Risk}\label{Siv}
Our goal in this subsection is to find a sharp asymptotic formula for
the Value at Risk ($\operatorname{VaR}_{\alpha}$) of the portfolio described
in Section~\ref{rm.sec}. The price $X_t$ at time $t$ for this portfolio
is defined by
(\ref{EBS}). We study the case where the confidence level $\alpha$
tends to one, and restrict ourselves to the portfolios with only
positive weights. The case of portfolios with both positive and
negative weights can be handled similarly.

For a portfolio, the value at risk $ \operatorname{VaR}_{\alpha}$, $0<\alpha<
1$, over the time period $t> 0$ is defined as the smallest number $k$
such that the probability of a loss greater than $k$ over the time
interval $t$ is equal to $\alpha$. It is not hard to see that
\[
\operatorname{VaR}_{\alpha}=\inf\bigl\{k\dvt \mathbb{P}(X_t\le
X_0-k)=1-\alpha\bigr\}.
\]

The next theorem provides an asymptotic formula for the function
$\alpha\mapsto \operatorname{VaR}_{\alpha}$ as the confidence level $\alpha$ tends to one.

\begin{theorem}\label{Tvarf}
Suppose assumption $(\mathcal A)$ holds for the covariance matrix
$\cm$. Then the following asymptotic formula is valid:
%
\begin{eqnarray}
\mathrm{VaR}_{\alpha} &=& X_0-\exp \biggl\{-\sqrt{
\frac
{2t}{\bar{A}_1+\cdots+\bar{A}_{\bar{n}}}\log\frac{1}{1-\alpha}}\nonumber \\
&&\label{Evarr}\hspace*{43pt}{}+\frac{\sum_{k=1}^{\bar{n}}\bar{A}_k
(\log({\bar{A}_1+\cdots+\bar{A}_{\bar{n}}})/{\bar{A}_k}+\mu_{k,t} )}{\bar{A}_1+\cdots+\bar{A}_{\bar{n}}} \biggr\}
\\
&&{}\times \biggl(1+\mathrm{O} \biggl(\frac{\log\log{1}/({1-\alpha})}{\sqrt{\log{1}/({1-\alpha})}} \biggr) \biggr)\nonumber
\end{eqnarray}
as $\alpha\rightarrow1$.
\end{theorem}

\begin{pf}
Let us fix $t> 0$, and denote by $F^{-1}_t$ the generalized inverse
function of the function
$F_t(x)=\mathbb{P}(X_t\le x)$. Then we have
%
\begin{equation}\label{Evar81}
\operatorname{VaR}_{\alpha}=X_0-F^{-1}_t(1-
\alpha).
\end{equation}
Therefore, in order to characterize the asymptotic behavior of the function
$\alpha\mapsto \operatorname{VaR}_{\alpha}$ as $\alpha\rightarrow1$, it
suffices to find an asymptotic formula
for the function $y\mapsto F^{-1}_t(y)$ as $y\rightarrow0$.

We will first study the asymptotics near zero of the inverse function
$F^{-1}$ of any function $F$, having the following form:
%
\begin{equation}\label{Evar1}
F(x)=c_1 \biggl(\frac{1}{x} \biggr)^{c_2} \biggl(
\log\frac{1}{x} \biggr)^{c_3} \exp \biggl\{-c_4
\log^2\frac{1}{x} \biggr\} \biggl(1+\mathrm{O} \biggl( \biggl(\log
\frac
{1}{x} \biggr)^{-1} \biggr) \biggr)
\end{equation}
as $x\rightarrow0$. It is assumed in (\ref{Evar1}) that the constants
satisfy the following conditions:
$c_1> 0$, $c_2\in\mathbb{R}$, $c_3\in\mathbb{R}$, and $c_4> 0$. We also
assume the continuity and the invertibility of the function $F$ near zero.

\begin{lemma}\label{Lfinal1}
Under the previous restrictions, the following asymptotic formula holds
as $y\rightarrow0$:
%
\begin{equation}\label{Evar10}
F^{-1}(y)=\exp\bigl\{-\sqrt{\phi(y)}\bigr\} \biggl(1+\mathrm{O} \biggl( \biggl(
\log\frac
{1}{y} \biggr)^{-1} \biggr) \biggr),
\end{equation}
where
\begin{eqnarray*}
\phi(y) &= & \frac{1}{c_4}\log\frac{1}{y}+c_2c_4^{-{3}/{2}}
\sqrt{\log \frac{1}{y}} +\frac{c_3}{2c_4}\log\log\frac{1}{y}
 \biggl(\frac{c_2^2}{2c_4^2}+\frac{c_3}{2c_4}\log\frac{1}{c_4}+
\frac
{1}{c_4}\log c_1 \biggr) \\
&&{}+\frac{c_2c_3}{4}c_4^{-{3}/{2}}
\frac{\log\log{1}/{y}}{\sqrt{\log{1}/{y}}}.
\end{eqnarray*}
\end{lemma}

Using the Taylor formula for the function $u\mapsto\sqrt{1+u}$ with two
terms, we obtain a simpler formula from (\ref{Evar10}).

\begin{corollary}\label{Cfinal2}
The following asymptotic formula holds:
%
\begin{equation}\label{Evar11}
F^{-1}(y)=\exp \biggl\{-\sqrt{\frac{1}{c_4}\log\frac{1}{y}}-\frac
{c_2}{2c_4} \biggr\} \biggl(1+\mathrm{O} \biggl(
\frac{\log\log{1}/{y}}{\sqrt{\log{1}/{y}}} \biggr) \biggr)
\end{equation}
as $y\rightarrow0$.
\end{corollary}

Note that formula (\ref{Evar11}) uses only the constants $c_2$ and $c_4$.

\begin{pf*}{Proof of Lemma~\ref{Lfinal1}}
Let $y> 0$, and let $F(u_y)=y$. Then $F^{-1}(y)=u_y$. Next, using (\ref
{Evar1}), we obtain
%
\begin{equation}\label{Evar2}
\log\frac{1}{y}=c_4\log^2\frac{1}{u_y}-c_2
\log\frac{1}{u_y}-c_3\log\log \frac{1}{u_y} -\log
c_1+\mathrm{O} \biggl( \biggl(\log\frac{1}{u_y} \biggr)^{-1}
\biggr)
\end{equation}
as $y\rightarrow0$. The previous formula implies the following
two-sided estimate:
\[
a_1\sqrt{\log\frac{1}{y}}\le\log\frac{1}{u_y}\le
a_2\sqrt{\log\frac
{1}{y}}, \qquad 0< y< y_0,
\]
for some constants $a_1> 0$ and $a_2> 0$.

Put $z_y=\log^2\frac{1}{u_y}$. Then formula (\ref{Evar2}) gives
%
\begin{equation} \label{Evar3}
z_y=\frac{1}{c_4}\log\frac{1}{y}+\frac{c_2}{c_4}
\sqrt{z_y}+\frac
{c_3}{2c_4}\log z_y +
\frac{1}{c_4}\log c_1+\mathrm{O} \biggl( \biggl(\log\frac{1}{y}
\biggr)^{-{1}/{2}} \biggr)
\end{equation}
as $y\rightarrow0$.

Our next goal is to use iterations in formula (\ref{Evar3}). We will
replace any occurrence of
$z_y$ on the right-hand side of (\ref{Evar3}) by the whole expression
on the right-hand side of
(\ref{Evar3}). The following simple formulas will be needed in the
sequel: $\log(1+s)=\mathrm{O}(s)$
and $\sqrt{1+s}=1+\frac{1}{2}s+\mathrm{O}(s^2)$ as $s\rightarrow0$. Let us put
%
\begin{equation}\label{Eh}
h=\frac{c_3}{2c_4}\log z_y+\frac{1}{4}\log
c_1+\mathrm{O} \biggl( \biggl(\log\frac
{1}{y} \biggr)^{-{1}/{2}}
\biggr),
\end{equation}
where the $\mathrm{O}$-term is the same as in formula (\ref{Evar3}).
We have
%
\begin{eqnarray}
\log z_y &=& \log \biggl(\frac{1}{c_4}\log\frac{1}{y}+
\frac{c_2}{c_4}\sqrt {z_y}+h \biggr)
\nonumber
\\[-8pt]
\label{Evar4}\\[-8pt]
\nonumber
&=& \log\frac{1}{c_4}+\log\log\frac{1}{y}+\mathrm{O} \biggl( \biggl(\log
\frac
{1}{y} \biggr)^{-{1}/{2}} \biggr)
\end{eqnarray}
as $y\rightarrow0$. Moreover,
\begin{eqnarray*}
\sqrt{z_y}&=& \frac{1}{\sqrt{c_4}}\sqrt{\log\frac{1}{y}}
\sqrt{1+\frac
{c_2\sqrt{z_y}}{\log{1}/{y}} +\frac{c_4h}{\log{1}/{y}}}
\\
&=&\frac{1}{\sqrt{c_4}}\sqrt{\log\frac{1}{y}}+\frac{c_2}{2\sqrt
{c_4}}
\frac{\sqrt{z_y}}{
\sqrt{\log{1}/{y}}}+\frac{\sqrt{c_4}h}{2\sqrt{\log{1}/{y}}} +\mathrm{O} \biggl( \biggl(\log\frac{1}{y}
\biggr)^{-{1}/{2}} \biggr)
\end{eqnarray*}
as $y\rightarrow0$. Next, we iterate again, and obtain
%
\begin{eqnarray}
\sqrt{z_y} &=& \frac{1}{\sqrt{c_4}}\sqrt{\log\frac{1}{y}}
\nonumber
\\
&&{}+\frac{c_2}{2\sqrt{c_4}\sqrt{\log{1}/{y}}}\nonumber
\\[-8pt]
\label{Evar5}\\[-8pt]
\nonumber
 &&\quad {}\times\biggl[\frac{1}{\sqrt{c_4}}\sqrt{\log\frac{1}{y}}+
\frac{c_2\sqrt
{z_y}}{2\sqrt{c_4}\sqrt{\log{1}/{y}}} +\frac{\sqrt{c_4}h}{2\sqrt{\log{1}/{y}}}+\mathrm{O} \biggl( \biggl(\log\frac
{1}{y}
\biggr)^{-{1}/{2}} \biggr) \biggr]
\\
&&{}+\frac{\sqrt{c_4}\hat{h}}{2\sqrt{\log{1}/{y}}}+\mathrm{O} \biggl( \biggl(\log \frac{1}{y}
\biggr)^{-{1}/{2}} \biggr)\nonumber
\end{eqnarray}
as $y\rightarrow0$. In (\ref{Evar5}), the symbol $\hat{h}$ stands for
the result of substituting
the expression for $\log z_y$ given in (\ref{Evar4}) into formula (\ref
{Eh}). It is not hard to see that
%
\begin{equation}\label{Evar6}
\frac{\sqrt{c_4}\hat{h}}{2\sqrt{\log{1}/{y}}}=\frac{c_3}{4\sqrt
{c_4}}\frac{\log\log{1}/{y}}{
\sqrt{\log{1}/{y}}}+\mathrm{O} \biggl( \biggl(\log
\frac{1}{y} \biggr)^{-{1}/{2}} \biggr)
\end{equation}
as $y\rightarrow0$. Now, taking into account (\ref{Evar5}) and (\ref
{Evar6}), we get
%
\begin{equation}\label{Evar7}
\sqrt{z_y}=\frac{1}{\sqrt{c_4}}\sqrt{\log\frac{1}{y}}+
\frac{c_2}{2c_4}+ \frac{c_3}{4c_4}\frac{\log\log{1}/{y}}{\sqrt{\log{1}/{y}}} +\mathrm{O} \biggl( \biggl(\log
\frac{1}{y} \biggr)^{-{1}/{2}} \biggr)
\end{equation}
as $y\rightarrow0$.
Finally, we can estimate $z_y$. Using (\ref{Evar3}), (\ref{Evar4}),
and (\ref{Evar7}),
we see that
%
\begin{eqnarray}
z_y&= &\frac{1}{c_4}\log\frac{1}{y}+c_2c_4^{-{3}/{2}}
\sqrt{\log\frac{1}{y}} +\frac{c_3}{2c_4}\log\log\frac{1}{y}
 \biggl(\frac{c_2^2}{2c_4^2}+\frac{c_3}{2c_4}\log\frac{1}{c_4}+
\frac
{1}{c_4}\log c_1 \biggr)\nonumber\\
&&\label{Evar8}{}+\frac{c_2c_3}{4}c_4^{-{3}/{2}}
\frac{\log\log{1}/{y}}{\sqrt{\log{1}/{y}}} +\mathrm{O} \biggl( \biggl(\log\frac{1}{y}
\biggr)^{-{1}/{2}} \biggr)
\\
&=& \phi(y)+\mathrm{O} \biggl( \biggl(\log\frac{1}{y} \biggr)^{-{1}/{2}}
\biggr) \nonumber
\end{eqnarray}
as $y\rightarrow0$.

We will next find an asymptotic formula for the function $F^{-1}$. We
will use the formula
%
\begin{equation}\label{Etake}
F^{-1}(y)=u_y=\exp\{-\sqrt{z_y}\}
\end{equation}
and the following simple lemma.

\begin{lemma}\label{Lsimp}
Let $z_y=\phi(y)+\mathrm{O}(\psi(y))$ as $y\rightarrow0$, where the functions
$\phi$ and $\psi$ are positive
and such that $\phi(y)\rightarrow\infty$ and $\frac{\psi(y)}{\sqrt{\phi
(y)}}\rightarrow0$ as
$y\rightarrow0$. Then
\[
u_y=\exp\bigl\{-\sqrt{\phi(y)}\bigr\} \biggl(1+\mathrm{O} \biggl(
\frac{\psi(y)}{\sqrt{\phi
(y)}} \biggr) \biggr)
\]
as $y\rightarrow0$.
\end{lemma}

\begin{pf*}{Proof of Lemma~\ref{Lsimp}}
We have
\begin{eqnarray*}
-\sqrt{z_y}&=&-\sqrt{\phi(y)+\mathrm{O}\bigl(\psi(y)\bigr)}=-\sqrt{\phi(y)}
\sqrt{1+\mathrm{O} \biggl(\frac{\psi(y)}{\phi(y)} \biggr)}
\\
&=& -\sqrt{\phi(y)} \biggl(1+\mathrm{O} \biggl(\frac{\psi(y)}{\phi(y)} \biggr) \biggr),
\end{eqnarray*}
and hence,
\begin{eqnarray*}
u_y&=&\exp\{-\sqrt{z_y}\}=\exp\bigl\{-\sqrt{\phi(y)}
\bigr\}\exp \biggl\{\mathrm{O} \biggl(\frac
{\psi(y)}{\sqrt{\phi(y)}} \biggr) \biggr\}
\\
&=&\exp\bigl\{-\sqrt{\phi(y)}\bigr\} \biggl(1+\mathrm{O} \biggl(\frac{\psi(y)}{\sqrt{\phi
(y)}}
\biggr) \biggr).
\end{eqnarray*}
\upqed\end{pf*}

Now, taking into account (\ref{Etake}), (\ref{Evar8}), and Lemma~\ref
{Lsimp} with the function $\phi$ such as in (\ref{Evar8})
and the function $\psi$ given by $\psi(y)= (\log\frac{1}{y}
)^{-{1}/{2}}$, we establish Lemma~\ref{Lfinal1}.
\end{pf*}

We are now ready to complete the proof of Theorem~\ref{Tvarf}.
Applying Theorem~\ref{gen.thm} to the random variable
$X_t$ given by (\ref{EBS}), we see that condition (\ref{Evar1})
holds. Note that
\[
c_2=-\frac{1}{t}\sum_{k=1}^{\bar{n}}
\bar{A}_k \biggl(\log\frac{\bar
{A}_1+\cdots+\bar{A}_{\bar{n}}}{\bar{A}_k} +\bar{\mu}_{m,t}
\biggr)\quad  \mbox{and}\quad  c_4=\frac{\bar{A}_1+\cdots+\bar{A}_{\bar{n}}}{2t}.
\]
Now, it is easy to see that (\ref{Evar81}) and Corollary~\ref{Cfinal2}
imply formula (\ref{Evarr}).
\end{pf}

\section*{Acknowledgements}
We would like to thank the two anonymous referees and the
editor in chief Eric Moulines, whose insightful remarks lead to a
substantial improvement of the paper.






\printhistory
\end{document}